\documentclass{article}
\usepackage{arxiv}
\usepackage[utf8]{inputenc} 
\usepackage[T1]{fontenc}    
\usepackage{hyperref}       
\usepackage{url}            
\usepackage{booktabs}       
\usepackage{amsfonts}       
\usepackage{nicefrac}       
\usepackage{microtype}      
\usepackage{lipsum}

\usepackage{algorithm}
\usepackage{algorithmic}
\usepackage{enumerate}
\usepackage{multirow}

\newtheorem{theorem}{Theorem}
\newtheorem{definition}{Definition}
\newtheorem{proposition}{Proposition}
\newtheorem{lemma}{Lemma}

\usepackage{amssymb,amsmath,latexsym}
\usepackage{graphicx}
\usepackage{caption}
\usepackage{subcaption}
\title{Sparse recovery based on the generalized error function}

\author{
  Zhiyong Zhou 
  \thanks{This work is supported by the Zhejiang Provincial Natural Science Foundation of China under Grant No.
  LQ21A010003.} \\
  Department of Statistics,
  Zhejiang University City College,
   Hangzhou, 310015, China \\
  \texttt{zhiyongzhou@zucc.edu.cn} \\
}

\begin{document}
\maketitle

\begin{abstract}
In this paper, we propose a novel sparse recovery method based on the generalized error function. The penalty function introduced involves both the shape and the scale parameters, making it very flexible. The theoretical analysis results in terms of the null space property, the spherical section property and the restricted invertibility factor are established for both constrained and unconstrained models. The practical algorithms via both the iteratively reweighted $\ell_1$ and the difference of convex functions algorithms are presented. Numerical experiments are conducted to illustrate the improvement provided by the proposed approach in various scenarios. Its practical application in magnetic resonance imaging (MRI) reconstruction is studied as well.
\end{abstract}


\section{Introduction}
High dimensionality is a basic feature of big data, that is, the number of features measured can be very large and are often considerable larger than the number of observations. To overcome the "curse of dimensionality" and reduce the redundancy, the vector of parameters to be estimated or the signal to be recovered is often assumed to be sparse (i.e., it has only a few nonzero entries) either by itself or after a proper transformation. How to exploit the sparsity to help estimating the underlying vector of parameters or recovering the unknown signal of interest, namely sparse recovery, has become a core research issue and gained immense popularity in the past decades \cite{hastie2015statistical}. 

Generally, sparse recovery aims to estimate an unknown sparse $\mathbf{x}\in\mathbb{R}^N$ from few noisy linear observations or measurements $\mathbf{y}=A\mathbf{x}+\boldsymbol{\varepsilon}\in\mathbb{R}^{m}$ where $A\in\mathbb{R}^{m\times N}$ with $m\ll N$ is the design or measurement matrix, and $\lVert \boldsymbol{\varepsilon}\rVert_2\leq \eta$ is the vector of noise. It arises in many scientific research fields including high-dimensional linear regression \cite{buhlmann2011statistics} and compressed sensing \cite{donoho2006compressed,eldar2012compressed,foucart2013mathematical}. Naturally, this sparse $\mathbf{x}$ can be recovered by solving a constrained $\ell_0$-minimization problem \begin{align}
\min\limits_{\mathbf{z}\in\mathbb{R}^N}\, \lVert \mathbf{z}\rVert_0\quad \text{subject to \quad $\lVert A\mathbf{z}-\mathbf{y}\rVert_2\leq \eta $}. \label{l0-mini}
\end{align}
or the unconstrained $\ell_0$-penalized least squares problem $\min_{\mathbf{z}\in\mathbb{R}^N} \frac{1}{2}\lVert A\mathbf{z}-\mathbf{y}\rVert_2^2+\lambda\lVert \mathbf{z}\rVert_0$ \cite{mazumder2011sparsenet}, where $\lambda>0$ is a tuning parameter. However, due to the nonsmoothness and nonconvexity of the $\ell_0$-norm, these are combinatorial problems which are known to be related to the selection of best subset and are computationally NP-hard to solve \cite{natarajan1995sparse}. Instead, a widely used solver is the following constrained $\ell_1$-minimization problem (also called Basis Pursuit Denoising) \cite{donoho2006compressed}:
\begin{align}
\min\limits_{\mathbf{z}\in\mathbb{R}^N}\, \lVert \mathbf{z}\rVert_1\quad \text{subject to \quad $\lVert A\mathbf{z}-\mathbf{y}\rVert_2\leq \eta $},
\end{align}
or the well-known Lasso  $\min_{\mathbf{z}\in\mathbb{R}^N} \frac{1}{2}\lVert A\mathbf{z}-\mathbf{y}\rVert_2^2+\lambda\lVert \mathbf{z}\rVert_1$ \cite{tibshirani1996regression}. The $\ell_1$-minimization acts as a convex relaxation of $\ell_0$-minimization. Although it enjoys attractive theoretical properties and has achieved great success in practice, it is biased and suboptimal. The Lasso does not have the oracle property \cite{fan2001variable}, since the $\ell_1$-norm is just a loose approximation of the $\ell_0$-norm.

To remedy this problem, many nonconvex sparse recovery methods have been employed to better approximate the $\ell_0$-norm and enhance sparsity. They include $\ell_p$ ($0<p<1$) \cite{chartrand2008iteratively,chartrand2008restricted,foucart2009sparsest}, Capped-L1 \cite{zhang2010analysis},  transformed $\ell_1$ (TL1) \cite{zhang2018minimization}, smooth clipped absolute deviation (SCAD) \cite{fan2001variable}, minimax concave penalty (MCP) \cite{zhang2010nearly}, exponential-type penalty (ETP) \cite{gao2011feasible,malek2016successive}, error function (ERF) method \cite{guo2020novel}, $\ell_1-\ell_2$ \cite{lou2018fast,yin2015minimization}, $\ell_r^r-\alpha\ell_1^r (\alpha\in[0,1],r\in(0,1])$ \cite{zhou2019new}, $\ell_1/\ell_2$ \cite{rahimi2018scale,wang2020accelerated}, and $q$-ratio sparsity minimization \cite{zhou2020minimization}, among others. These parameterized nonconvex methods result in the difficulties of theoretical analysis and computational algorithms due to the nonconvexity of the penalty functions, but do outperform the convex $\ell_1$-minimization in various scenarios. For example, it has been reported that $\ell_p$ gives superior results for incoherent measurement matrices, while $\ell_1-\ell_2$, $\ell_1/\ell_2$ and $q$-ratio sparsity minimization are better choices for highly coherent measurement matrices. Meanwhile, TL1 is a robust choice no matter whether the measurement matrix is coherent or not. 

The resulting nonconvex sparse recovery methods have the following general constrained form: \begin{align}
\min\limits_{\mathbf{z}\in\mathbb{R}^N}\, R_{\theta}(\mathbf{z}), \quad \text{subject to \quad $\lVert A\mathbf{z}-\mathbf{y}\rVert_2\leq \eta $}, \label{general}
\end{align}
where $R_{\theta}(\cdot):\mathbb{R}^N\rightarrow\mathbb{R}_{+}:=[0,\infty)$ denotes a nonconvex penalty or regularization function with an approximation parameter $\theta$ (it can be a vector of parameters), or a formulation of its corresponding unconstrained version $ \min_{\mathbf{z}\in\mathbb{R}^N} \frac{1}{2}\lVert A\mathbf{z}-\mathbf{y}\rVert_2^2+\lambda R_{\theta}(\mathbf{z})$ with $\lambda>0$ being the tuning parameter. Several nonconvex methods and their corresponding penalty functions are shown in Table \ref{table:1}. Basically, in practice a separable and concave on $\mathbb{R}_{+}^N$ penalty function is desired in order to facilitate the theoretical analysis and solving algorithms. 

\begin{table}[h!]
\centering
 \begin{tabular}{||c| c| c ||} 
 \hline
 Method &  Penalty &  Property \\ [0.5ex] 
 \hline\hline
 $\ell_p$ & $\lVert \mathbf{z}\rVert_p^p, p\in(0,1)$ & \multirow{13}{*}{ separable, concave on $\mathbb{R}_{+}^N$}\\ 
 \cline{1-2}
 Capped-L1 & $\sum\limits_{j=1}^N \min\{|z_j|,\gamma\}, \gamma>0$ & \\
 \cline{1-2}
 TL1 &  $\sum\limits_{j=1}^N \frac{(a+1)|z_j|}{a+|z_j|}, a>0$ & \\
  \cline{1-2}
 MCP & $\sum\limits_{j=1}^N \int_{0}^{|z_j|}\lambda (1-\frac{t}{\gamma\lambda})_{+} dt, \lambda>0, \gamma>0$ &  \\
  \cline{1-2}
 SCAD & $\sum\limits_{j=1}^N \int_{0}^{|z_j|} \lambda\left(1\wedge (1-\frac{t-\lambda}{\lambda(\gamma-1)})_{+} \right)dt, \lambda>0, \gamma>1$ &  \\ 
  \cline{1-2}
 ETP & $\sum\limits_{j=1}^N 1-e^{-|z_j|/\tau}, \tau>0$ &\\
  \cline{1-2}
 ERF & $\sum\limits_{j=1}^N \int_{0}^{|z_j|} e^{-t^2/\sigma^2} dt, \sigma>0$ & \\
 \hline
 $\ell_1-\ell_2$ & $\lVert \mathbf{z}\rVert_1-\lVert \mathbf{z}\rVert_2$ & \multirow{4}{*} {nonseparable, nonconcave on $\mathbb{R}_{+}^N$} \\
 \cline{1-2}
 $\ell_r^r-\alpha\ell_1^r$ & $\lVert \mathbf{z}\rVert_r^r-\alpha\lVert \mathbf{z}\rVert_1^r, \alpha\in [0,1],r\in(0,1]$ & \\
  \cline{1-2}
 $\ell_1/\ell_2$ & $\lVert \mathbf{z}\rVert_1/\lVert \mathbf{z}\rVert_2$ &  \\
  \cline{1-2}
 $q$-ratio sparsity & $s_q(\mathbf{z})=(\lVert \mathbf{z}\rVert_1/\lVert \mathbf{z}\rVert_q)^{\frac{q}{q-1}}, q\in[0,\infty]$ &  \\ [1ex] 
 \hline
\end{tabular} 
\caption{Examples of nonconvex methods together with their corresponding penalties and element properties.}
\label{table:1}
\end{table}

For the nonconvex methods mentioned above, there are a large number of theoretical recovery results for their global or local optimal solutions (please refer to the corresponding literature for a specific method). Moreover, some unified recovery analysis results have also been obtained for these nonconvex sparse recovery methods. For instance, \cite{tran2019class} established a theoretical recovery guarantee through unified null space properties. A general theoretical framework was presented in \cite{zhang2012general} based on regularity conditions. It shows that under appropriate conditions, the global solution of concave regularization leads to desirable recovery performance and the global solution corresponds to the unique sparse local solution.

From an algorithmic point of view, many optimization algorithms have been used for solving the nonconvex problem (\ref{general}) and its unconstrained version. Among them, \cite{fan2001variable} proposed a local quadratic approximation (LQA) for the SCAD method. A local linear approximation (LLA) was suggested in \cite{zou2008one}, which is a special case of the iteratively reweighted $\ell_1$ (IRL1) algorithm  \cite{candes2008enhancing,zhao2012reweighted}. The iterative reweighted algorithms were generalized in \cite{ochs2015iteratively} to tackle a problem class of the sum of a convex function and a nondecreasing function applied to another convex function. In \cite{zhang2010analysis}, a general framework of multi-stage convex relaxation was given, which iteratively refines the convex relaxation to achieve better solutions for nonconvex optimization problems. It includes previous approaches including LQA, LLA and the concave-convex procedure (CCCP) \cite{yuille2003concave} as special cases. Under appropriate conditions, the author shows that the local solution obtained by the multi-stage convex relaxation method for nonconvex regularization achieves better recovery performance than the $\ell_1$-regularization. In addition, the proximal gradient method \cite{parikh2014proximal} has also been extended to address the nonconvex problem (\ref{general}). The ingredient of this method hinges on computing the proximal operator of the penalty function. It is well-known that the proximal operators of the $\ell_0$-norm and the $\ell_1$-norm are the hard-thresholding operator and the soft-thresholding operator, respectively. Most of nonconvex penalty functions also have proximal operators with closed-form. But for some penalty functions such as the ERF, there is no proximal operator with closed-form such that the corresponding proximal gradient method is inefficient in practice. Another general methodology that can be adopted to solve (\ref{general}) and its unconstrained version is referred to the Difference of Convex functions Algorithms (DCA) \cite{tao1997convex,tao1998dc}. DCA reformulate the nonconvex penalty function $R_{\theta}(\cdot)$ as a difference of two convex functions and then solve it based on the DC programming, please see \cite{le2015dc} for a comprehensive view of DCA approaches for sparse optimization.

\subsection{Contribution}
In this paper, we introduce a very flexible nonconvex sparse recovery method with the ETP and ERF approaches as special cases. We recognize that the generalized error function can be used to construct a smooth concave penalty function and provide a unified framework. The main contributions of the present paper are four-folds. 
\begin{enumerate}[(i)]
    \item We propose a general sparse recovery approach based on the generalized error function. The properties and the proximal operator of the proposed penalty are studied.
    \item We establish a rigorous recovery analysis for the proposed approach in terms of the null space property, the spherical section property and the restricted invertibility factor for both constrained and unconstrained models. 
    \item We present efficient algorithms to solve the proposed method via both IRL1 and DCA, and carry out various numerical experiments to illustrate their superior performance.
    \item We apply the proposed approach to magnetic resonance imaging (MRI) reconstruction, which improves the performance of state-of-the-art methods.
\end{enumerate}

\subsection{Related Work}
The first study of the exponential-type penalty (ETP) was in the support vector machines problem \cite{bradley1998feature}. Then this penalty together with its modifications was introduced into sparse regression \cite{gao2011feasible}, low-rank matrix recovery \cite{malek2013recovery,malek2014iterative} and sparse signal recovery \cite{malek2016successive}. Very recently, an approach based on the error function (ERF) has been proposed in \cite{guo2020novel} to promote sparsity for signal recovery. In comparison to these mostly related works, our main distinctions are summarized as follows:
\begin{itemize}
\item We have provided a fairly general framework of sparse recovery based on the generalized error function, making the existing ETP and ERF methods as its special cases.
\item We have systematically studied the theoretical recovery analysis results and practical solving algorithms of our proposed method, while the existing literatures only contain part of our results. For instance, we have generalized the spherical section property based analysis in \cite{malek2016successive} from $\ell_2$-bound to $\ell_q$-bound for any $q\in(1,\infty]$. We have discussed the solving algorithms based on both IRL1 and DCA, while only IRL1 was considered in \cite{guo2020novel} for the ERF method. 
\item To the best of our knowledge, we are the first to adopt this kind of exponential-type approaches to MRI reconstruction.
\end{itemize}

\subsection{Terminologies and Outline}
We use the following terminologies throughout the paper. We denote vectors by lower bold case letters and matrices by upper case letters. Vectors are columns by defaults. $x_j$ denotes the $j$-th component of the vector $\mathbf{x}$. $\mathrm{sign}(x_j)$ is the sign function which is $x_j/|x_j|$ if $x_j\neq 0$ and zero otherwise. We introduce the notations $[N]$ for the set $\{1,2,\cdots,N\}$ and $|S|$ for the cardinality of a set $S$ in $[N]$. We write $S^c$ for the complement $[N]\setminus S$ of $S$. For any vector $\mathbf{x}\in\mathbb{R}^N$, its $\ell_q$ "norm" is defined as $\lVert \mathbf{x}\rVert_q=(\sum_{j=1}^N |x_j|^q)^{1/q}$ for any $q\in(0,\infty)$, with the extension $\lVert \mathbf{x}\rVert_0=\sum_{j=1}^N 1_{\{x_j\neq 0\}}$ and $\lVert \mathbf{x}\rVert_\infty=\max_{1\leq j\leq N} |x_j|$. For a vector $\mathbf{x}\in\mathbb{R}^N$, we say $\mathbf{x}$ is $k$-sparse if at most $k$ of its entries are nonzero, i.e., if $\lVert \mathbf{x}\rVert_0\leq k$. $\mathbf{x}_S$ denotes the vector which coincides with $\mathbf{x}$ on the indices in $S$ and is extended to zero outside $S$. $\mathrm{supp}(\mathbf{x})$ denotes the support set of $\mathbf{x}$. For any matrix $A\in\mathbb{R}^{m\times N}$, the kernel of $A$ is defined as $\mathrm{Ker}(A)=\{\mathbf{x}\in\mathbb{R}^N|A\mathbf{x}=\mathbf{0}\}$ and $A_{S}$ denotes the submatrix with columns restricted to the index set $S$. We use $\mathrm{tr}(A)$ to denote the trace of a square matrix $A$.

An outline of this paper is as follows. In Section 2, we introduce the proposed nonconvex sparse recovery method. We study the theoretical recovery analysis in Section 3. In Section 4, both the iteratively reweighted $\ell_1$ and the difference of convex functions algorithms are used to solve the nonconvex regularization problem. The proposed method adapted to the gradient for imaging applications is also discussed. In Section 5, a series of numerical experiments are conducted. Conclusions and some future works are presented in Section 6.

\section{Proposed Method}

\begin{figure}[htbp]
	\centering
	\includegraphics[width=\textwidth,height=0.4\textheight]{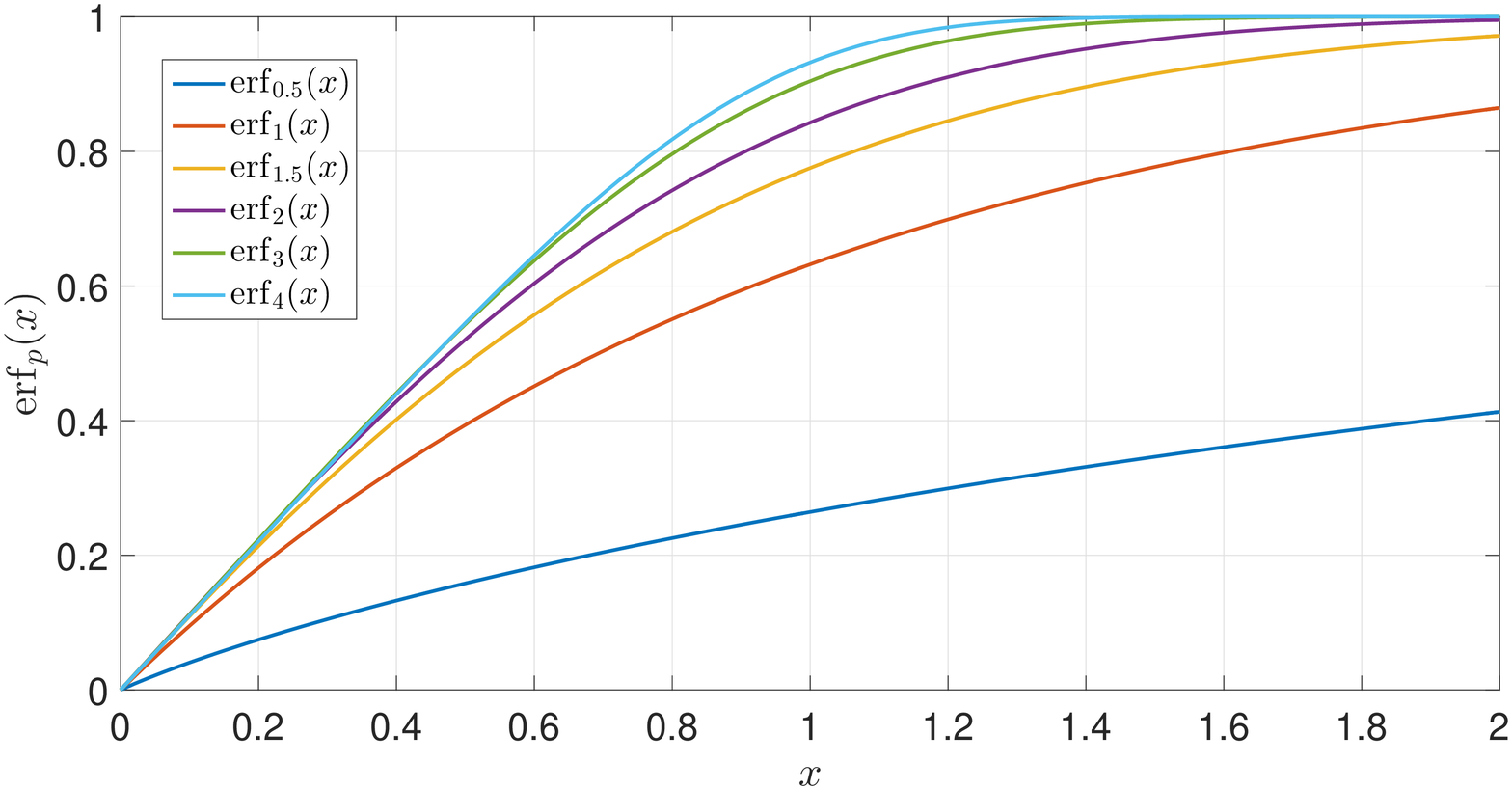}
	\caption{Plots for the generalized error functions $\mathrm{erf}_p(\cdot)$ with $p=0.5,1,1.5,2,3,4$.} \label{gerf}
\end{figure}

We start with the definition of the generalized error function. For any $p>0$, the generalized error function \cite{yin2018some,wang2019nonconvex} is defined as 
\begin{align}
    \mathrm{erf}_p(x)=\frac{p}{\Gamma(1/p)}\int_{0}^{x}{e^{-t^p}dt}.
\end{align}
Then, we have $\mathrm{erf}_1(x)=1-e^{-x}$ and $\mathrm{erf}_2(x)=\mathrm{erf}(x)=\frac{2}{\sqrt{\pi}}\int_{0}^{x}{e^{-t^2}dt}$ (note that this is the standard error function). In addition, it is easy to verify that for any given $p>0$, the generalized error function $\mathrm{erf}_p(x)$ is concave on $\mathbb{R}_{+}$, strictly increasing with respect to $x\in\mathbb{R}_{+}$, $\mathrm{erf}_p(0)=0$ and $\lim_{x\rightarrow\infty}\mathrm{erf}_p(x)=1$. These properties can be observed from the function plots for the $\mathrm{erf}_p(\cdot)$ in $\mathbb{R}_{+}$ with different choices of $p$ displayed in Figure \ref{gerf}.

Based on this generalized error function (GERF), for any $\mathbf{x}\in\mathbb{R}^N$ and $p>0$, we can use the following penalty function or regularization term for sparse recovery \begin{align}
    J_{p,\sigma}(\mathbf{x})=\sum\limits_{j=1}^N \Phi_{p,\sigma}(|x_j|),
\end{align}
where $\Phi_{p,\sigma}(x)=\int_{0}^{x}{e^{-(\tau/\sigma)^p}d\tau}=\sigma\int_{0}^{x/\sigma}{e^{-t^p}dt}=\frac{\sigma\Gamma(1/p)}{p}\mathrm{erf}_p(\frac{x}{\sigma})$. Here $p>0$ is the shape parameter, while $\sigma>0$ is the scale parameter. It is obvious that when $p\rightarrow 0$, then $\Phi_{p,\sigma}(x)\rightarrow \Phi_{0,\sigma}(x)=x$ and hence $J_{p,\sigma}(\mathbf{x})$ reduces to the $\ell_1$-norm $\lVert \mathbf{x}\rVert_1$. The parameters $p$ and $\sigma$ control the degree of concavity of the penalty function $J_{p,\sigma}(\cdot)$. We refer our proposed sparse recovery approach based on this penalty function $J_{p,\sigma}(\cdot)$ as "GERF" throughout this paper. It is very flexible and includes some existing sparsity-promoting models as special cases. In particular, the exponential-type penalty (ETP) introduced in \cite{gao2011feasible,malek2016successive} is the special case of $p=1$, and the sparse recovery method based on the error function (referred as ERF) in \cite{guo2020novel} corresponds to the special case of $p=2$.

In fact, regardless of some constant, our GERF penalty $J_{p,\sigma}(\cdot)$ is an interpolation of the $\ell_0$ and $\ell_1$ penalty. It gains smoothness over the $\ell_0$ penalty and hence allows more computational options. Moreover, it improves the variable selection accuracy and achieves the oracle property by reducing the bias of the $\ell_1$-minimization. The ability to choose suitable values of both the shape parameter $p$ and the scale parameter $\sigma$ enables us to fully exploit its sparsity-promoting power and to achieve better recovery performance.

\subsection{Properties of the Penalty}

In what follows, we provide some useful properties for the functions $\Phi_{p,\sigma}(\cdot):\mathbb{R}\rightarrow\mathbb{R}$ and $J_{p,\sigma}(\cdot):\mathbb{R}^N\rightarrow\mathbb{R}_{+}$.
\begin{itemize}
    \item  The derivative of the function $\Phi_{p,\sigma}(\cdot)$ goes to \begin{align*}
        \frac{d}{dx} \Phi_{p,\sigma}(x)=e^{-(x/\sigma)^p}.
    \end{align*}
    \item $\Phi_{p,\sigma}(x)$ is concave and non-decreasing in $x\in\mathbb{R}_{+}$.
    \item According to the following double inequality given in \cite{alzer1997some}: \begin{align*}
        (1-e^{-bx^p})^{1/p}<\frac{1}{\Gamma(1+1/p)}\int_{0}^x e^{-t^p} dt<(1-e^{-ax^p})^{1/p}, 
    \end{align*}
    where $x>0$, $p\neq 1$ is a positive real number and \[
  \begin{cases}
    a=1 \quad \text{and}\quad b= [\Gamma(1+1/p)]^{-p}, \quad 0<p<1 \\
    a=[\Gamma(1+1/p)]^{-p}\quad \text{and}\quad b=1, \quad p>1
  \end{cases}
    \]
    we are able to get the corresponding double inequality for $\Phi_{p,\sigma}(\cdot)$.
    \item $J_{p,\sigma}(-\mathbf{x})=J_{p,\sigma}(\mathbf{x})$ for any $\mathbf{x}\in\mathbb{R}^N$ and $J_{p,\sigma}(\mathbf{0})=0$.
    \item $J_{p,\sigma}(\cdot)$ is concave on $\mathbb{R}_{+}^N$, that is for any $\mathbf{x},\mathbf{y}\in\mathbb{R}_{+}^N$ and any $t\in[0,1]$, we have \begin{align*}
       t J_{p,\sigma}(\mathbf{x})+(1-t)J_{p,\sigma}(\mathbf{y})\leq  J_{p,\sigma}(t\mathbf{x}+(1-t)\mathbf{y})
    \end{align*}
    \item $J_{p,\sigma}(\cdot)$ is sub-additive so that for any $\mathbf{x},\mathbf{y}\in\mathbb{R}^N$,  \begin{align*}
        J_{p,\sigma}(\mathbf{x}+\mathbf{y})\leq  J_{p,\sigma}(\mathbf{x})+ J_{p,\sigma}(\mathbf{y}).
    \end{align*}
    And the additive property holds when $\mathbf{x},\mathbf{y}\in\mathbb{R}^N$ have disjoint supports.
\end{itemize}

Next, we present a proposition that characterizes the asymptotic behaviors of $J_{p,\sigma}(\cdot)$ with respect to $p$ and $\sigma$.

\begin{proposition}
For any nonzero $\mathbf{x}\in\mathbb{R}^N$, $p>0$ and $\sigma>0$, \begin{enumerate}[(a)]
    \item It is trivial that  $J_{p,\sigma}(\mathbf{x})\rightarrow \lVert \mathbf{x}\rVert_1$, as $p\rightarrow 0^{+}$.
    \item $J_{p,\sigma}(\mathbf{x})\rightarrow \lVert \mathbf{x}\rVert_1$, as $\sigma\rightarrow +\infty$. 
    \item $J_{p,\sigma}(\mathbf{x})/\sigma\rightarrow \frac{\Gamma(1/p)}{p}\lVert \mathbf{x}\rVert_0$, as $\sigma\rightarrow 0^{+}$. Moreover, if $p\rightarrow \infty$, then we have $J_{p,\sigma}(\mathbf{x})/\sigma\rightarrow \lVert \mathbf{x}\rVert_0$ since $\lim\limits_{p\rightarrow \infty}\frac{\Gamma(1/p)}{p}=\lim\limits_{p\rightarrow \infty}\Gamma(\frac{1}{p}+1)=\Gamma(1)=1$.
\end{enumerate}
\end{proposition}

{\bf Proof.} With regard to (b), when $x\neq 0$, it holds that \begin{align*}
\lim\limits_{\sigma\rightarrow +\infty} \frac{\Phi_{p,\sigma}(x)}{x}=\lim\limits_{t\rightarrow 0}\frac{\int_{0}^{t}{e^{-x^p}dx}}{t}=1,
\end{align*}
where we let $t=x/\sigma$. When $x=0$, it is obvious that $\Phi_{p,\sigma}(x)=0=x$.\\

As for (c), we have $\Phi_{p,\sigma}(0)=0$ and $\forall\, x\neq 0$, \begin{align*}
    \lim\limits_{\sigma\rightarrow 0^{+}} \frac{\Phi_{p,\sigma}(x)}{\sigma}=\int_{0}^{\infty}{e^{-t^p}dt}=\frac{1}{p}\int_{0}^{\infty}{x^{1/p-1}e^{-x}dx}=\frac{\Gamma(1/p)}{p}.
\end{align*}

\begin{figure}[htbp]
	\centering
	\includegraphics[width=\textwidth,height=0.4\textheight]{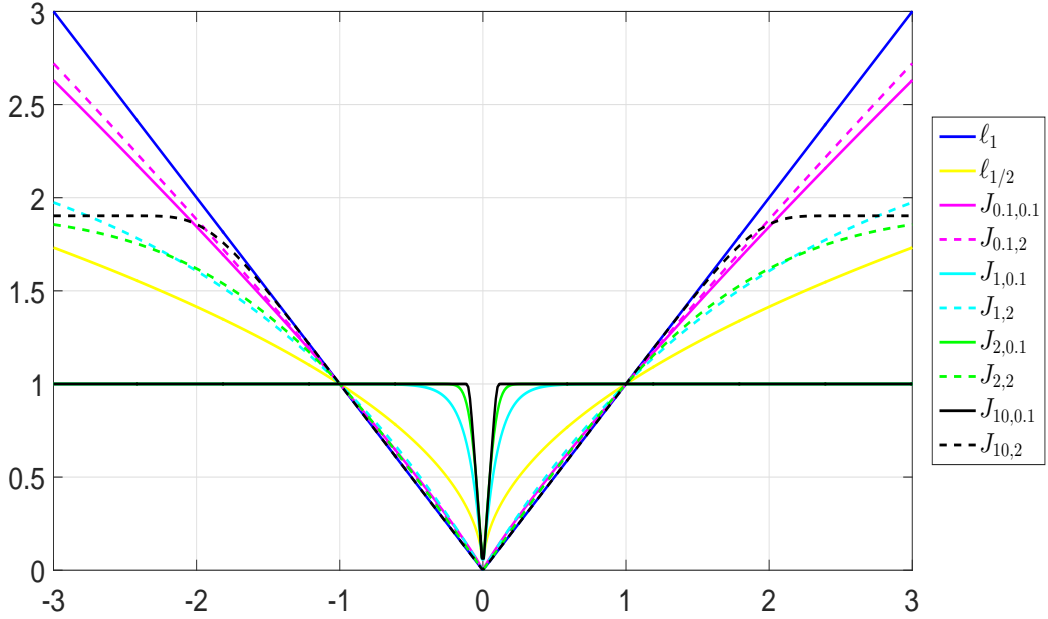}
	\caption{The objective functions of the proposed GERF penalty $J_{p,\sigma}(\cdot)$ in $\mathbb{R}$ with various values of $p$ and $\sigma$, compared to the $\ell_1$ and $\ell_{1/2}$ penalties.} \label{ObjectiveFunction}
\end{figure}

The objective functions of the GERF penalty $J_{p,\sigma}(\cdot)$ in $\mathbb{R}$  with various values of $p$ and $\sigma$ are displayed in Figure \ref{ObjectiveFunction}. All the functions are scaled to attain the point $(1,1)$ for a better comparison. As shown, the proposed GERF model gives a good approximation to the $\ell_0$ penalty for a small value of $\sigma$, and approaches the $\ell_1$ penalty for a small value of $p$ or a large value of $\sigma$.

\begin{figure}[htbp]
	\centering
	\includegraphics[width=\textwidth,height=0.35\textheight]{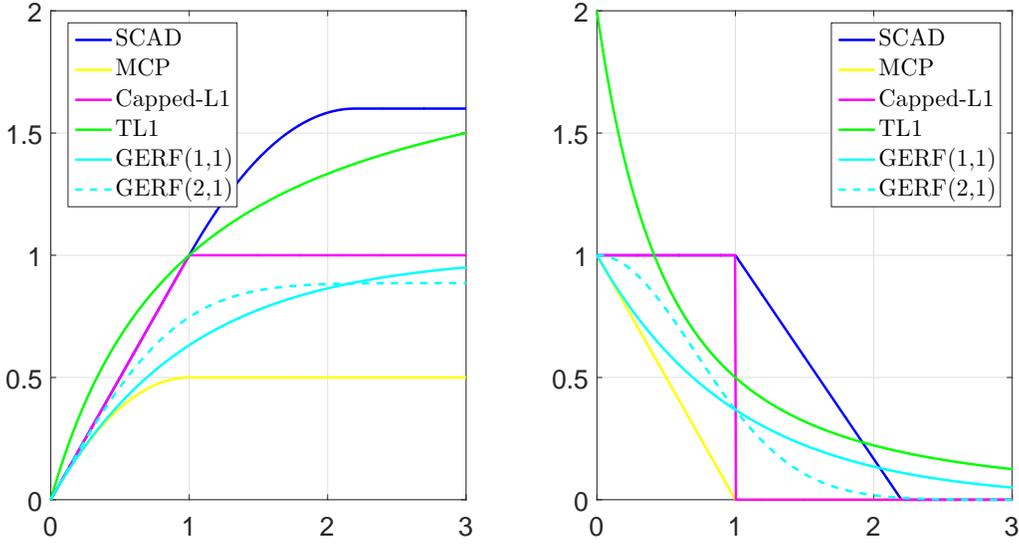}
	\caption{The curves of the SCAD ($\gamma=2.2$), the MCP ($\gamma=1$), the Capped-L1 ($\gamma=1$), the TL1 ($a=1$), the GERF with $p=1,\sigma=1$ and $p=2,\sigma=1$ (Left Panel), and their derivatives (Right Panel) in $\mathbb{R}_{+}$.} \label{penalty}
	\end{figure}

Moreover, we compare the GERF penalty functions ($p=1,\sigma=1$ and $p=2,\sigma=1$) and their derivatives with other popular penalty functions in $\mathbb{R}_{+}$. From Figure \ref{penalty}, it can be seen that the GERF penalty functions have very similar patterns as other widely used penalty functions. Specifically, all the penalties displayed are concave, nondecreasing on $\mathbb{R}_{+}$, and are folded or near-folded (i.e., the functions become constant for large variable) \cite{liu2017folded}. Their derivatives are all nonincreasing, finite at $0$, and equal to or close to zero when the variable of function is large, which satisfy the criterion for an unbiased penalty function that its derivative has a horizontal asymptote at zero. What stand out for our proposed GERF penalty are its smoothness and flexibility.

\subsection{Proximal operator}
In this subsection, we discuss the proximal operator for the GERF. For a function $\psi:\mathbb{R}^N\rightarrow \mathbb{R}$ and a parameter $\mu>0$, the proximal operator $\mathrm{prox}_{\mu\psi}$ is defined as \begin{align}
\mathrm{prox}_{\mu\psi}(\mathbf{x})=\mathop{\arg\min}_{\mathbf{u}\in\mathbb{R}^N}\left\{\frac{1}{2\mu}\lVert \mathbf{u}-\mathbf{x}\rVert_2^2+\psi(\mathbf{u})\right\}.
\end{align}
As we know, the proximal operator for the $\ell_1$-norm is the soft-thresholding operator $S_{\mu}(\mathbf{x})=\mathrm{sign}(\mathbf{x})(|\mathbf{x}|-\mu)_{+}$, while the proximal operator for the $\ell_0$-norm is the hard-thresholding operator $H_{\mu}(\mathbf{x})=\mathbf{x}1_{|\mathbf{x}|>\mu}$. With regard to the GERF model, the proximal operator for the penalty function $J_{p,\sigma}(\cdot)$ fulfills the optimality condition as follows: \begin{align*}
    \mathbf{x}\in \mu\exp(-(|\mathbf{u}|/\sigma)^p)\partial |\mathbf{u}|+\mathbf{u}.
\end{align*}

\begin{figure}[htbp]
	\centering
	\includegraphics[width=\textwidth,height=0.4\textheight]{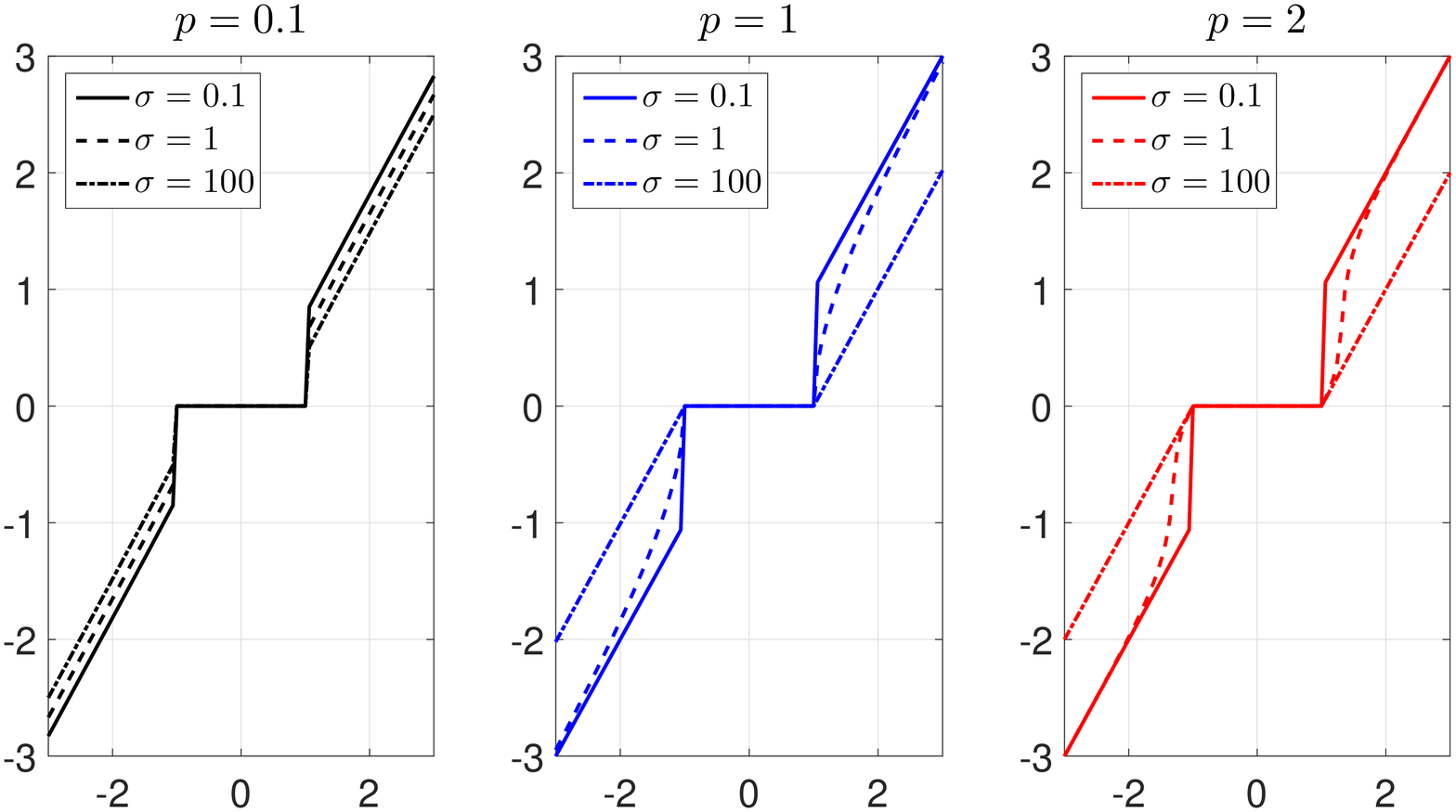}
	\caption{Proximal operators of the GERF while varying $p$ and $\sigma$ with $\mu=1$.} \label{proximal_operator}
\end{figure}

Then, it is easy to observe that we have $u_j=0$ when $|x_j|\leq \mu$, otherwise the minimizer is given as the solution of the following nonlinear equation: \begin{align*}
    x_j=\mu\exp(-(|u_j|/\sigma)^p)\mathrm{sign}(x_j)+u_j.
\end{align*}
In general, the solution has no closed-form but can be numerically found via the Newton's method. It is worth mentioning that for the special case of $p=1$ (i.e., for the ETP method), a solution with closed-form as well as the corresponding proximal operator can be obtained by using the \emph{Lambert W} function \cite{corless1996lambertw}, please see the details in \cite{malek2016successive}. 
Here we present the proximal operators for the GERF with $p=0.1,1,2$ and $\sigma=0.1,1,100$ in Figure \ref{proximal_operator} when $\mu$ is fixed to $1$. As we can see, for a large $\sigma=100$, the proximal operators for the GERF are close to the soft-thresholding operator, whereas they are similar to the hard-thresholding operator for a small $\sigma=0.1$. These findings are largely consistent with earlier statements in Proposition 1.

\section{Recovery Analysis}
In this section, we set out to establish the theoretical recovery analysis for the proposed GERF approach. Both constrained and unconstrained models are considered.
\subsection{Null Space Property}
We begin with the null space property based analysis for the following noiseless constrained GERF minimization problem:
 \begin{align}
\min\limits_{\mathbf{z}\in\mathbb{R}^N} J_{p,\sigma}(\mathbf{z})\quad  \text{subject to} \quad A\mathbf{z}=A\mathbf{x}. \label{noiseless_gerf}
\end{align}

It is well-known that an exact sparse recovery can be guaranteed for the noiseless $\ell_1$-minimization if and only if the measurement matrix $A$ satisfies a null space property (NSP) \cite{cohen2009compressed}: \begin{align}
\lVert \mathbf{v}_S\rVert_1< \lVert \mathbf{v}_{S^c}\rVert_1,  \quad \forall\,\, \mathbf{v}\in \mathrm{Ker}(A)\setminus \{\mathbf{0}\}.
\end{align}

Similarly, a generalized version of NSP guarantees an exact sparse recovery for our proposed GERF approach.

\begin{definition}
Given $p>0$ and $\sigma>0$, we say a matrix $A\in\mathbb{R}^{m\times N}$ satisfies a generalized null space property (gNSP) relative to $J_{p,\sigma}(\cdot)$ and $S\subset [N]$ if \begin{align}
    J_{p,\sigma}(\mathbf{v}_S)< J_{p,\sigma}(\mathbf{v}_{S^c}) \quad \text{for all $\mathbf{v}\in\mathrm{Ker}(A)\setminus \{\mathbf{0}\}$}.
\end{align}
It satisfies the gNSP of order $s$ relative to $J_{p,\sigma}(\cdot)$ if it satisfies the gNSP relative to $J_{p,\sigma}(\cdot)$ and any $S\subset [N]$ with $|S|\leq s$
\end{definition}

Then, based on the sub-additive property of the penalty function $J_{p,\sigma}(\cdot)$, it is straightforward to obtain the following theorem, which provides a sufficient and necessary condition for an exact sparse recovery via the GERF approach. Its proof follows almost the same as that given in Section 3.3 of \cite{guo2020novel} and so is omitted.

\begin{theorem}
For any given measurement matrix $A\in\mathbb{R}^{m\times N}$, $p>0$ and $\sigma>0$, every $s$-sparse vector $\mathbf{x}\in\mathbb{R}^N$ is the unique solution of the problem (\ref{noiseless_gerf})
if and only if $A$ satisfies the gNSP of order $s$ relative to $J_{p,\sigma}(\cdot)$.
\end{theorem}

We note that since the penalty function $J_{p,\sigma}(\cdot)$ is separable and symmetric on $\mathbb{R}^N$, and concave on $\mathbb{R}_{+}^N$, the following result holds as an immediate consequence of Proposition 4.6 in \cite{tran2019class}, which verifies that our proposed GERF approach is superior to the $\ell_1$-minimization in sparse uniform recovery.

\begin{theorem}
Let $N>1$, $s\in\mathbb{N}$ with $1\leq s<N/2$. Then we have \begin{align*}
    \left\{\mathbf{v}\in\mathbb{R}^N:\lVert \mathbf{v}_S\rVert_1<\lVert \mathbf{v}_{S^c}\rVert_1,\forall\, S\subset [N]\, \text{with}\, |S|\leq s\right\}\subset  \left\{\mathbf{v}\in\mathbb{R}^N: J_{p,\sigma}(\mathbf{v}_S)<
    J_{p,\sigma}(\mathbf{v}_{S^c}),\forall\, S\subset [N]\, \text{with}\, |S|\leq s\right\}.
\end{align*}
Consequently, the gNSP relative to $J_{p,\sigma}(\cdot)$ is less demanding than the NSP for $\ell_1$-minimization. 
\end{theorem}

\subsection{Spherical Section Property}
Regarding the noiseless model (\ref{noiseless_gerf}), in this subsection we further show that by choosing a proper $\sigma$, one can obtain a solution arbitrarily close to the unique solution of $\ell_0$-minimization via (\ref{noiseless_gerf}) when $p$ is fixed.

\begin{definition}
For any $q\in(1,\infty]$, we say the measurement matrix $A$ possesses the $\Delta_q$-spherical section property if $\Delta_q(A)\leq (\lVert \mathbf{v}\rVert_1/\lVert \mathbf{v}\rVert_q)^{\frac{q}{q-1}}$ holds for all $\mathbf{v}\in \mathrm{Ker}(A)\setminus \{\mathbf{0}\}$. Hence, we have \begin{align}
    \Delta_q(A)=\inf_{\mathbf{v}\in \mathrm{Ker}(A)\setminus \{\mathbf{0}\}} \left(\frac{\lVert \mathbf{v}\rVert_1}{\lVert \mathbf{v}\rVert_q}\right)^{\frac{q}{q-1}}
\end{align}
\end{definition}

\noindent
{\bf Remark.} The $\Delta_q$-spherical section property introduced here is a direct extension of $\Delta$-spherical section property studied in \cite{malek2016successive,zhang2013theory} from $q=2$ to any $q\in(1,\infty]$. Based on this newly defined spherical section property, the results established in this subsection generalize the corresponding results in \cite{malek2016successive,zhang2013theory}.

\begin{proposition}
Suppose $A\in\mathbb{R}^{m\times N}$ has the $\Delta_q$-spherical property for some $q\in(1,\infty]$. If $\lVert \mathbf{x}\rVert_0<2^{\frac{q}{1-q}}\Delta_q(A)$, then $\mathbf{x}$ is the unique solution of the following noiseless $\ell_0$-minimization problem: \begin{align}
    \min\limits_{\mathbf{z}\in\mathbb{R}^N}\,\lVert \mathbf{z}\rVert_{0}\quad  \text{subject to} \quad A\mathbf{z}=A\mathbf{x}. \label{noiseless_l0}
\end{align}
\end{proposition}

{\bf Proof.}  If there exists a vector $\tilde{\mathbf{x}}$ ($\tilde{\mathbf{x}}\neq \mathbf{x}$) such that $A\tilde{\mathbf{x}}=A\mathbf{x}$ and $\lVert \tilde{\mathbf{x}} \rVert_0\leq \lVert \mathbf{x}\rVert_0$. Then $\mathbf{v}=\tilde{\mathbf{x}}-\mathbf{x}\in \mathrm{Ker}(A)\setminus \{\mathbf{0}\}$. Therefore, if $\lVert \mathbf{x}\rVert_0<2^{\frac{q}{1-q}}\Delta_q(A)$ holds for some $q\in(1,\infty]$, then we have \begin{align*}
    \left(\frac{\lVert \mathbf{v}\rVert_1}{\lVert \mathbf{v}\rVert_q}\right)^{\frac{q}{q-1}}\geq \Delta_q(A)>2^{\frac{q}{q-1}}\lVert \mathbf{x}\rVert_0.
\end{align*}
However, $\lVert \mathbf{v}\rVert_0=\lVert\tilde{\mathbf{x}}-\mathbf{x}\rVert_0\leq \lVert \tilde{\mathbf{x}} \rVert_0   
+ \lVert \mathbf{x}\rVert_0\leq 2\lVert \mathbf{x}\rVert_0$. Moreover, it holds that $ \left(\frac{\lVert \mathbf{v}\rVert_1}{\lVert \mathbf{v}\rVert_q}\right)^{\frac{q}{q-1}}\leq \lVert \mathbf{v}\rVert_0$. As a result, we  arise in a contradiction that \begin{align*}
    2\lVert \mathbf{x}\rVert_0\geq \lVert \mathbf{v}\rVert_0\geq \left(\frac{\lVert \mathbf{v}\rVert_1}{\lVert \mathbf{v}\rVert_q}\right)^{\frac{q}{q-1}}>2^{\frac{q}{q-1}}\lVert \mathbf{x}\rVert_0.
\end{align*}
Hence, $\mathbf{x}$ is the unique solution of the $\ell_0$-minimization problem (\ref{noiseless_l0}).

\begin{lemma}
Suppose $A\in\mathbb{R}^{m\times N}$ has the $\Delta_q$-spherical section property for some $q\in(1,\infty]$. Let $\mathbf{v}\in \mathrm{Ker}(A)\setminus \{\mathbf{0}\} $. Then for any subset $S\subset [N]$ such that $|S|+\Delta_q(A)>N$, we have \begin{align}
    \frac{\lVert\mathbf{v}_S\rVert_1}{\lVert \mathbf{v}\rVert_q}\geq \Delta_q(A)^{1-1/q}-(N-|S|)^{1-1/q}.\label{bound_lemma1}
\end{align}
Moreover, if $\mathbf{v}$ has at most $\lceil\Delta_q(A)-1\rceil$ entries with absolute values greater than $\alpha$, then \begin{align}
    \lVert \mathbf{v}\rVert_q \leq \frac{N\alpha}{\Delta_q(A)^{1-1/q}-(\lceil\Delta_q(A)-1\rceil)^{1-1/q}}, \label{bound_2}
\end{align}
where $\lceil c \rceil$ denotes the smallest integer greater than or equal to $c$.
\end{lemma}

{\bf Proof.} When $S=[N]$, this lemma holds trivially since $(\lVert\mathbf{v}_S\rVert_1/\lVert \mathbf{v}\rVert_q)^{q/(q-1)}\geq \Delta_q(A)$. Otherwise, when $|S|<N$, it holds that \begin{align*}
    \Delta_q(A)^{1-1/q}\leq \frac{ \lVert\mathbf{v}\rVert_1}{\lVert\mathbf{v}\rVert_q}=\frac{\lVert\mathbf{v_S}\rVert_1+\lVert\mathbf{v_{S^c}}\rVert_1}{\lVert\mathbf{v}\rVert_q}.
\end{align*}
Meanwhile, $\lVert\mathbf{v_{S^c}}\rVert_1\leq (N-|S|)^{1-1/q}\lVert\mathbf{v_{S^c}}\rVert_q\leq (N-|S|)^{1-1/q}\lVert\mathbf{v}\rVert_q $, which completes the proof of (\ref{bound_lemma1}).

To prove (\ref{bound_2}), we let $S$ be the indices of entries with absolute values not greater than $\alpha$, then $|S|\geq N-\lceil\Delta_q(A)-1\rceil$. By using (\ref{bound_lemma1}), we have 
\begin{align*}
    \lVert \mathbf{v}\rVert_q\leq \frac{\lVert \mathbf{v}_S\rVert_1}{\Delta_q(A)^{1-1/q}-(N-|S|)^{1-1/q}}\leq \frac{N\alpha}{\Delta_q(A)^{1-1/q}-(\lceil\Delta_q(A)-1\rceil)^{1-1/q}}.
\end{align*}
This proof is completed. 

\begin{lemma}
Assume $A\in\mathbb{R}^{m\times N}$ has the $\Delta_q$-spherical section property for some $q\in(1,\infty]$ and $\lVert \mathbf{x}\rVert_0=k< 2^{\frac{q}{1-q}}\Delta_q(A)$. Then for any $\hat{\mathbf{x}}\in \{\mathbf{z}\in\mathbb{R}^N:A\mathbf{z}=A\mathbf{x}\}$ and $\alpha_{p,\sigma}=\Phi_{p,\sigma}^{-1}\left(\frac{\sigma\Gamma(1/p)}{p}(1-\frac{1}{N})\right)=\sigma\mathrm{erf}_p^{-1}(1-\frac{1}{N})$, we have \begin{align}
    \lVert \hat{\mathbf{x}}-\mathbf{x}\rVert_q\leq \frac{N\alpha_{p,\sigma}}{\Delta_q(A)^{1-1/q}-(\lceil\Delta_q(A)-1\rceil)^{1-1/q}},
\end{align}
whenever $J_{p,\sigma}(\hat{\mathbf{x}})\leq \frac{\sigma\Gamma(1/p)}{p}(\lceil \Delta_q(A)-1\rceil-k)$.
\end{lemma}

{\bf Proof.} Assume $T$ denotes the indices of entries with absolute values greater than $\alpha_{p,\sigma}$ in $\hat{\mathbf{x}}$. Note that \begin{align*}
    J_{p,\sigma}(\hat{\mathbf{x}})&=J_{p,\sigma}(\hat{\mathbf{x}}_T)+J_{p,\sigma}(\hat{\mathbf{x}}_{T^c}) \\
    &=\sum\limits_{j\in T} \Phi_{p,\sigma}(|\hat{x}_j|)+\sum\limits_{j\in T^c} \Phi_{p,\sigma}(|\hat{x}_j|)\\
    &> \frac{\sigma\Gamma(1/p)}{p} |T|(1-1/N)>\frac{\sigma\Gamma(1/p)}{p}(|T|-1).
\end{align*}
Meanwhile, $J_{p,\sigma}(\hat{\mathbf{x}})\leq \frac{\sigma\Gamma(1/p)}{p}(\lceil \Delta_q(A)-1\rceil-k)$, which yields $|T|<(\lceil \Delta_q(A)-1\rceil-k)+1$ such that $|T|\leq \lceil \Delta_q(A)-1\rceil-k$. That is to say, there are at most $\lceil \Delta_q(A)-1\rceil-k$ entries of $\hat{\mathbf{x}}$ that have absolute values greater than $\alpha_{p,\sigma}$. Furthermore, since $\lVert \mathbf{x}\rVert_0=k$, so at most $k+\lceil \Delta_q(A)-1\rceil-k=\lceil \Delta_q(A)-1\rceil$ entries of $\hat{\mathbf{x}}-\mathbf{x}$ has absolute values that are greater than $\alpha_{p,\sigma}$. Then the fact the $\hat{\mathbf{x}}-\mathbf{x}\in \mathrm{Ker}(A)\setminus \{\mathbf{0}\}$ and Lemma 1 complete the proof of this lemma.

\begin{lemma}
If $\hat{\mathbf{x}}$ is the solution of the problem (\ref{noiseless_gerf}) and the conditions in Lemma 2 hold, then it satisfies \begin{align*}
    J_{p,\sigma}(\hat{\mathbf{x}})\leq \frac{\sigma\Gamma(1/p)}{p}(\lceil \Delta_q(A)-1\rceil-k).
\end{align*}
\end{lemma}

{\bf Proof.} Notice that \begin{align*}
    J_{p,\sigma}(\hat{\mathbf{x}})\leq  J_{p,\sigma}(\mathbf{x})=\sum\limits_{j=1}^N \Phi_{p,\sigma}(|x_j|)\leq \frac{\sigma\Gamma(1/p)}{p} k\leq \frac{\sigma\Gamma(1/p)}{p} (2^{\frac{q}{q-1}}k-k).
\end{align*}
In addition, $2^{\frac{q}{q-1}}k<\Delta_q(A)\leq  \lceil \Delta_q(A)-1\rceil+1$ yields that $2^{\frac{q}{q-1}}k\leq \lceil \Delta_q(A)-1\rceil$. Therefore, the proof is completed.

Then, as a result of Lemma 2 and Lemma 3, the following theorem holds immediately. 

\begin{theorem}
Let $\hat{\mathbf{x}}$ be the solution of the problem (\ref{noiseless_gerf}) and assume the conditions in Lemma 2 hold. Then $\hat{\mathbf{x}}$ approaches the unique solution of the noiseless $\ell_0$-minimization problem (\ref{noiseless_l0}) $\mathbf{x}$ as $\sigma\mathrm{erf}_{p}^{-1}(1-\frac{1}{N})$ approaches $0$. In particular, when the shape parameter $p$ is fixed, this happens when the scale parameter $\sigma$ approaches $0$.
\end{theorem}

\subsection{Restricted Invertibility Factor}
Next, we move on to discuss the recovery analysis result for the solution of the unconstrained or regularized version of GERF approach: \begin{align}
    \hat{\mathbf{x}}=\mathop{\arg\min}\limits_{\mathbf{z}\in\mathbb{R}^N}\,\frac{1}{2m}\lVert \mathbf{y}-A\mathbf{z}\rVert_2^2+\lambda J_{p,\sigma}(\mathbf{z}), \label{unconstrained_method}
\end{align}
where $\mathbf{y}=A\mathbf{x}+\boldsymbol{\varepsilon}\in\mathbb{R}^m$ with $\boldsymbol{\varepsilon}\sim N(0,\tilde{\sigma}^2 \mathrm{I}_{m})$. Following the general theory of concave regularization presented in \cite{zhang2012general}, our concave penalty has the form \begin{align}
    \lambda J_{p,\sigma}(\mathbf{z})=\sum\limits_{j=1}^N \rho(z_j;\lambda),
\end{align}
where $\rho(t;\lambda)=\lambda \int_{0}^{|t|}{e^{-(x/\sigma)^p}dx}$ is the scalar regularization function.

Then the threshold level of this penalty $\lambda^{*}:=\inf_{t>0}\{t/2+\rho(t;\lambda)/t\}$, which is a function of $\lambda$ that provides a normalization of $\lambda$, satisfies that $\lambda^{*}\leq \lim_{t\rightarrow 0^{+}} (\partial/\partial t)\rho(t;\lambda)=\lim_{t\rightarrow 0^{+}} \lambda e^{-(t/\sigma)^p}=\lambda$. In what follows, based on the tools of restricted invertibility factor and $\eta$ null consistency condition used in \cite{zhang2012general}, we are able to obtain an error bound for the GERF regularization method. For a rigorous proof of the main result, the reader is referred to \cite{zhang2012general} for detailed arguments.

\begin{definition} (\cite{zhang2012general})
For $q\geq 1$, $\xi>0$ and $S\subset[N]$, the restricted invertibility factor (RIF) is defined as \begin{align}
\mathrm{RIF}_{q}(\xi,S)=\inf\limits_{\{\mathbf{v}\in\mathbb{R}^N:J_{p,\sigma}(\mathbf{v}_{S^c})<\xi J_{p,\sigma}(\mathbf{v}_S)\}}\frac{|S|^{1/q}\lVert A^{T}A\mathbf{v}\rVert_{\infty}}{m\lVert \mathbf{v}\rVert_{q}}
\end{align}
\end{definition}

\begin{definition} (\cite{zhang2012general})
For a given $\eta\in(0,1]$, we say the regularization method presented in (\ref{unconstrained_method}) satisfies the $\eta$ null consistency condition ($\eta$-NC) if it holds that \begin{align}
    \min\limits_{\mathbf{z}\in\mathbb{R}^N} \{ \lVert \boldsymbol{\varepsilon}/\eta-A\mathbf{z}\rVert_2^2/(2m)+\lambda J_{p,\sigma}(\mathbf{z})\}=\lVert \boldsymbol{\varepsilon}/\eta\rVert_2^2/(2m).
\end{align}
\end{definition}

\begin{theorem}
Assume $S=\mathrm{supp}(\mathbf{x})$ and the $\eta$-NC condition holds with some $\eta\in(0,1)$. Then for all $q\geq 1$, we have \begin{align}
\lVert \hat{\mathbf{x}}-\mathbf{x}\rVert_q \leq \frac{\lambda(1+\eta) |S|^{1/q}}{\mathrm{RIF}_{q}(\frac{1+\eta}{1-\eta},S)}.
\end{align}
\end{theorem}

\section{Algorithms}
In this section, we present the algorithms to solve the proposed GERF method via both IRL1 and DCA, and study the GERF on the gradient for imaging applications. We consider the following unconstrained problem: \begin{align}
\min\limits_{\mathbf{x}\in\mathbb{R}^N} \frac{1}{2}\lVert \mathbf{y}-A\mathbf{x}\rVert_2^2+\lambda J_{p,\sigma}(\mathbf{x}), \label{unconstrained}
\end{align}
where $\lambda>0$ is the tuning parameter. Note that, in order to solve the constrained version of the GERF approach, we merely need to make some small modifications to the algorithms for solving the unconstrained model (\ref{unconstrained}). For the sake of brevity of the article, we will not repeat the discussion for the constrained model.

\subsection{IRL1}
As discussed in \cite{guo2020novel,ochs2015iteratively}, an iteratively reweighted $\ell_1$ (IRL1) algorithm can be used to solve the nonsmooth nonconvex problem (\ref{unconstrained}). At each iteration, this algorithm needs to solve a reweighted $\ell_1$ subproblem given as follows:
\begin{align}
   \mathbf{x}^{(k+1)}=\mathop{\arg\min}\limits_{\mathbf{x}\in\mathbb{R}^N} \frac{1}{2}\lVert \mathbf{y}-A\mathbf{x}\rVert_2^2+\lambda \sum\limits_{j=1}^N w_j^{(k)} |x_j| \label{reweighted_l1}
\end{align}
with weights $w_j^{(k)}=\exp(-|x_j^{(k)}/\sigma|^p)$. Then an alternating direction method of multipliers (ADMM) algorithm \cite{boyd2011distributed} can be adopted to solve this subproblem (\ref{reweighted_l1}). We introduce an auxiliary vector $\boldsymbol{\theta}$ and then consider the augmented Lagrangian
\begin{align}
\mathcal{L}_{\rho}(\mathbf{x},\boldsymbol{\theta};\boldsymbol{\beta})=\frac{1}{2}\lVert \mathbf{y}-A\mathbf{x}\rVert_2^2+\lambda \sum\limits_{j=1}^N w_j^k |x_j|+\frac{\rho}{2}\lVert \mathbf{x}-\boldsymbol{\theta}+\boldsymbol{\beta}\rVert_2^2,
\end{align}
where $\rho>0$ is a small fixed positive parameter and $\boldsymbol{\beta}$ is the dual variable. The
ADMM algorithm is based on minimizing the augmented Lagrangian $\mathcal{L}_{\rho}(\mathbf{x},\boldsymbol{\theta};\boldsymbol{\beta})$ successively over $\mathbf{x}$ and $\boldsymbol{\theta}$, and then applying a dual variable update to $\boldsymbol{\beta}$, which yields the following updates:\[
\begin{cases}
\mathbf{x}^{t+1}=\mathop{\arg\min}\limits_{\mathbf{x}}\mathcal{L}_{\rho}(\mathbf{x},\boldsymbol{\theta}^{t};\boldsymbol{\beta}^{t})=S_{\lambda\mathbf{w}^{(k)}/\rho}(\boldsymbol{\theta}^{t}-\boldsymbol{\beta}^{t}),\\
\boldsymbol{\theta}^{t+1}=\mathop{\arg\min}\limits_{\boldsymbol{\theta}}\mathcal{L}_{\rho}(\mathbf{x}^{t+1},\boldsymbol{\theta};\boldsymbol{\beta}^{t})=(A^{T}A+\rho\mathrm{I}_{N})^{-1}(A^{T}\mathbf{y}+\rho\mathbf{x}^{t+1}+\rho\boldsymbol{\beta}^{t}),\\
\boldsymbol{\beta}^{t+1}=\boldsymbol{\beta}^{t}+\mathbf{x}^{t+1}-\boldsymbol{\theta}^{t+1},
\end{cases}
\]
for iterations $t=0,1,2,\cdots$.

\subsection{DCA}

Another type of approaches that can be used to solve the nonsmooth nonconvex problem (\ref{unconstrained}) is the Difference of Convex functions Algorithms (DCA) \cite{tao1997convex,tao1998dc}, where we express the nonconvex penalty as a difference of two convex functions. With regard to our penalty $J_{p,\sigma}(\cdot)$, since \begin{align*}
    \Phi_{p,\sigma}(|x|)=\int_{0}^{|x|}{e^{-(\tau/\sigma)^p}d\tau}=|x|-(|x|-\int_{0}^{|x|}{e^{-(\tau/\sigma)^p}d\tau})=|x|-\int_{0}^{|x|}{(1-e^{-(\tau/\sigma)^p})d\tau},
\end{align*}
hence we have the following DC decomposition: \begin{align}
     J_{p,\sigma}(\mathbf{x})=\sum\limits_{j=1}^N \Phi_{p,\sigma}(|x_j|)=\lVert \mathbf{x}\rVert_1-\sum\limits_{j=1}^N \int_{0}^{|x_j|}{(1-e^{-(\tau/\sigma)^p})d\tau}.
\end{align}
Then the problem (\ref{unconstrained}) can be expressed as a DC program:\begin{align*}
    \min\limits_{\mathbf{x}\in\mathbb{R}^N} \frac{1}{2}\lVert \mathbf{y}-A\mathbf{x}\rVert_2^2+\lambda J_{p,\sigma}(\mathbf{x})&= \min\limits_{\mathbf{x}\in\mathbb{R}^N} \frac{1}{2}\lVert \mathbf{y}-A\mathbf{x}\rVert_2^2+\lambda\lVert \mathbf{x}\rVert_1-\lambda\sum\limits_{j=1}^N \int_{0}^{|x_j|}{(1-e^{-(\tau/\sigma)^p})d\tau} \\
    &=\min\limits_{\mathbf{x}\in\mathbb{R}^N} g(\mathbf{x})-h(\mathbf{x}),
\end{align*}
where $g(\mathbf{x})=\frac{1}{2}\lVert \mathbf{y}-A\mathbf{x}\rVert_2^2+\lambda\lVert \mathbf{x}\rVert_1$ and $h(\mathbf{x})=\lambda\sum\limits_{j=1}^N \int_{0}^{|x_j|}{(1-e^{-(\tau/\sigma)^p})d\tau}$.

The procedure of DCA for solving the problem (\ref{unconstrained}) is summarized in Algorithm \ref{algorithm_dca}.
\begin{algorithm}[!h]
	\caption{DCA for solving (\ref{unconstrained})} 
	\vspace*{0.5em}
	\begin{algorithmic}[1]
	    \item Input: Measurement matrix $A\in\mathbb{R}^{m\times N}$, Measurement vector $\mathbf{y}\in\mathbb{R}^m$, Parameters $p>0$, $\sigma>0$ and $\lambda>0$.
	    \item Output: The solution of the problem (\ref{unconstrained}).
	    
		\STATE {\textbf{Initialization}: Set $\mathbf{x}^{(0)}=\mathbf{0}$, $k=0$.}
		\STATE {\textbf{Iteration}: Repeat until the stopping rule is met (e.g., $k>Max$ or $\lVert \mathbf{x}^{(k+1)}-\mathbf{x}^{(k)}\rVert_2/\max\{\lVert \mathbf{x}^{(k)}\rVert_2,1\}<Tol$), \begin{align}
			\mathbf{x}^{(k+1)}=
			\mathop{\arg\min}\limits_{\mathbf{x}\in\mathbb{R}^N}  \frac{1}{2}\lVert \mathbf{y}-A\mathbf{x}\rVert_2^2+\lambda\lVert \mathbf{x}\rVert_1-\lambda\langle\mathbf{v}^{(k)},\mathbf{x}\rangle,
			\end{align}
			where $\mathbf{v}^{(k)}\in\mathbb{R}^N$ with entries $v_j^{(k)}=\mathrm{sign}(x_j^{(k)})\left(1-e^{-|x_j^{(k)}/\sigma|^p}\right)$ for $j\in[N]$.}
		\STATE {\textbf{Update iteration}: $k=k+1$.}
	\end{algorithmic} \label{algorithm_dca}
\end{algorithm} 

Note that, for each DCA iteration, it is required to solve the following $\ell_1$-regularized subproblem: \begin{align}
    \min\limits_{\mathbf{x}\in\mathbb{R}^N} \frac{1}{2}\lVert \mathbf{y}-A\mathbf{x}\rVert_2^2+\lambda\lVert \mathbf{x}\rVert_1+\langle\mathbf{v},\mathbf{x}\rangle.
\end{align}
Here we can also employ an ADMM algorithm to solve this subproblem by introducing an auxiliary vector $\boldsymbol{\theta}$ and then considering the augmented Lagrangian
\begin{align}
\mathcal{L'}_{\rho}(\mathbf{x},\boldsymbol{\theta};\boldsymbol{\beta})=\frac{1}{2}\lVert \mathbf{y}-A\mathbf{x}\rVert_2^2+\langle\mathbf{v},\mathbf{x}\rangle+\lambda\lVert \boldsymbol{\theta}\rVert_1+\boldsymbol{\beta}^{T}(\mathbf{x}-\boldsymbol{\theta})+
\frac{\rho}{2}\lVert \mathbf{x}-\boldsymbol{\theta}\rVert_2^2,
\end{align}
where $\rho>0$ is a small fixed positive parameter and $\boldsymbol{\beta}$ is the dual variable. The resulting ADMM algorithm consists of the following updates:\[
\begin{cases}
\mathbf{x}^{t+1}=\mathop{\arg\min}\limits_{\mathbf{x}}\mathcal{L'}_{\rho}(\mathbf{x},\boldsymbol{\theta}^{t};\boldsymbol{\beta}^{t})=(A^{T}A+\rho\mathrm{I}_{N})^{-1}(A^{T}\mathbf{y}-\mathbf{v}+\rho\boldsymbol{\theta}^{t}-\boldsymbol{\beta}^{t}),\\
\boldsymbol{\theta}^{t+1}=\mathop{\arg\min}\limits_{\boldsymbol{\theta}}\mathcal{L'}_{\rho}(\mathbf{x}^{t+1},\boldsymbol{\theta};\boldsymbol{\beta}^{t})=S_{\lambda/\rho}(\mathbf{x}^{t+1}+\boldsymbol{\beta}^{t}/\rho),\\
\boldsymbol{\beta}^{t+1}=\boldsymbol{\beta}^{t}+\rho(\mathbf{x}^{t+1}-\boldsymbol{\theta}^{t+1}),
\end{cases}
\]
for iterations $t=0,1,2,\cdots$.

\subsection{GERF on the gradient}
In this subsection, we discuss the imaging applications with a GERF minimization adapted to the gradient. Let $u\in\mathbb{R}^{m\times N}$ be the image of interest, $A$ be the measurement matrix and $f$ be the measurements. For the MRI reconstruction, $A=\mathcal{F}R$, where $\mathcal{F}$ denotes the Fourier transform and $R$ is the sampling mask in the frequency space. We denote $D_x, D_y$ as the horizontal and vertical partial derivative operators respectively, and $D=[D_x; D_y]$ is the gradient operator $\nabla$ in the discrete setting. Then, the reconstruction algorithm based on the GERF approach goes to \begin{align}
\min\limits_{u\in\mathbb{R}^{m\times N}}  J_{p,\sigma} (Du)\quad \text{subject to \quad $Au=f$}, \label{gerf_gradient}
\end{align}
where $Du=[D_x u;D_y u]$. Let $(u_{jx},u_{jy})$ be the gradient vector at pixel $j$, then $J_{p,\sigma}(Du)$ can be rewritten as \begin{align}
    J_{p,\sigma}(Du)=\sum\limits_{j} \left(\int_{0}^{|u_{jx}|} e^{-(\tau/\sigma)^p}d\tau +\int_{0}^{|u_{jy}|} e^{-(\tau/\sigma)^p}d\tau\right).
\end{align}

Then we are able to use the framework of DCA proposed previously to solve this problem. In each DCA iteration, the subproblem reduces to a total variation (TV) \cite{rudin1992nonlinear} type of minimization, which can be handled efficiently by using the split Bregman technique \cite{goldstein2009split}. In this paper, in order to solve (\ref{gerf_gradient}) we directly adopt a simple modification of the Algorithm 2 in \cite{lou2015weighted} (or Algorithm 3 in \cite{yin2015minimization}), where a combination of DCA and split Bregman technique is used for the $\ell_1-\ell_2$ minimization adapted to the gradient. We are only required to replace the approximation vector $q^{n}=(q_x^n, q_y^n)=(D_x u^n,D_y u^n)/\sqrt{|D_x u^n|^2+|D_y u^n|^2}$ in \cite{lou2015weighted} by a new $q^{n}=(q_x^n, q_y^n)=\left(\mathrm{sign}(D_x u^n).*(1-e^{-|D_x u^n/\sigma|^p}),\mathrm{sign}(D_y u^n).*(1-e^{-|D_y u^n/\sigma|^p})\right)$ for our GERF approach. For the sake of conciseness, we will not reproduce the procedure for this algorithm. The reader is referred to \cite{lou2015weighted} or \cite{yin2015minimization} for details.

\section{Numerical Experiments}
In this section, we conduct a series of numerical experiments to demonstrate the advantageous performance of our proposed GERF method. Throughout this section, the support of $k$-sparse signal is a random index set and its nonzero entries follow a standard normal distribution.

\subsection{A Test with IRL1 and DCA}
\begin{figure}[htbp]
	\centering
	\includegraphics[width=\textwidth,height=0.4\textheight]{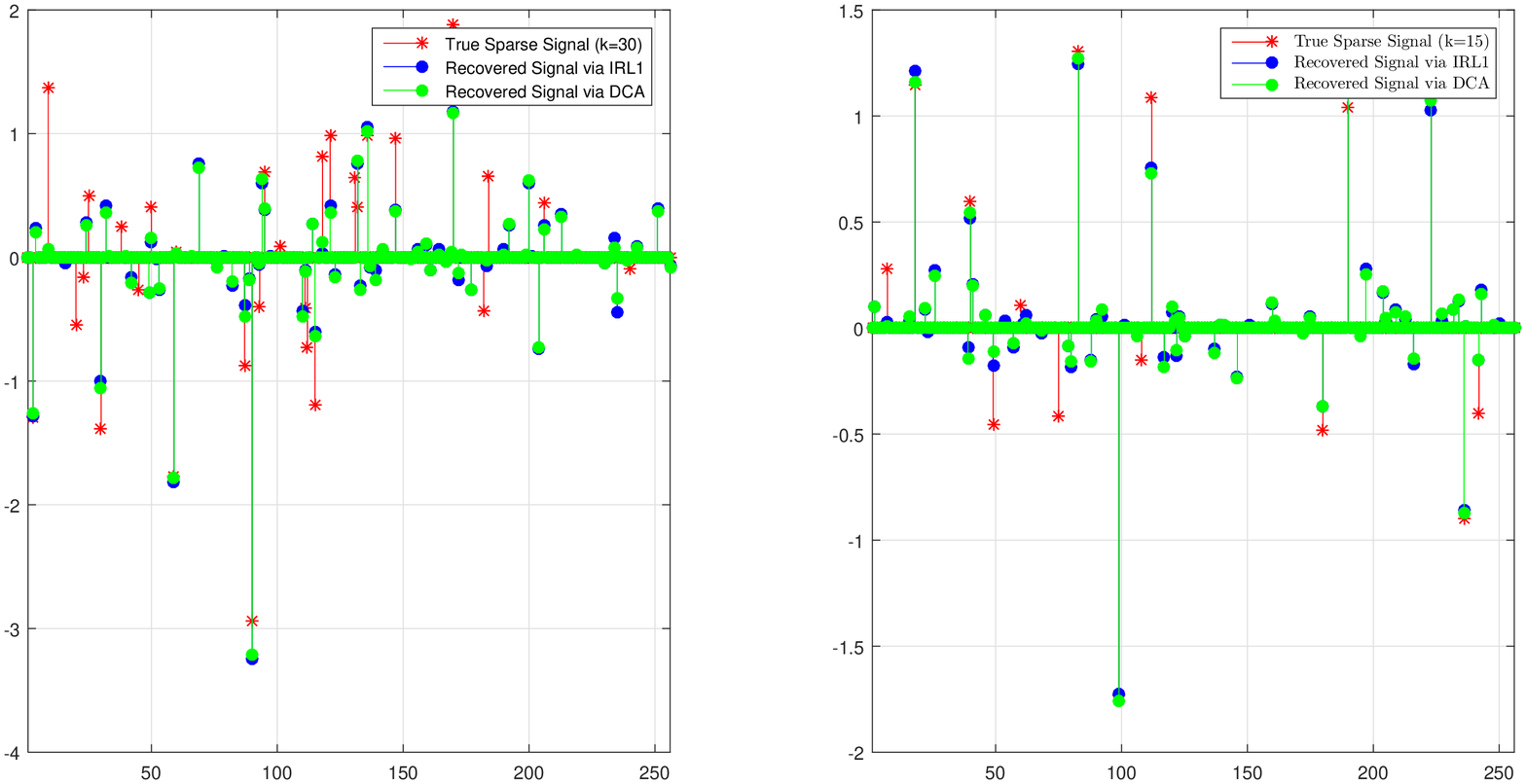}
	\caption{Reconstruction comparison between the IRL1 and the DCA for $p=2,\sigma=1$. Left panel (Noise Free) and Right panel (Noisy with $\tilde{\sigma}=0.1$).} \label{IRL1_DCA}
\end{figure}

Firstly, we do a numerical test to compare the reconstruction performance of the proposed GERF approach by using the IRL1 and DCA algorithms with $p=2,\sigma=1$ for a sparse signal $\mathbf{x}\in\mathbb{R}^{256}$ reconstruction with a Gaussian random measurement matrix $A\in\mathbb{R}^{64\times 256}$. Two cases are considered here, one is that the true signal $\mathbf{x}$ has a sparsity level of $30$ and the measurements are noise free, while the other case is that the true signal is $15$-sparse and the measurements are noisy with $\tilde{\sigma}=0.1$. We set $\lambda=10^{-5}$ and $\lambda=10^{-2}$ for the noise free and noisy cases, respectively. It can be observed from Figure \ref{IRL1_DCA} that in both cases the reconstructed signals via IRL1 and DCA methods are almost the same. The corresponding computing time for these two algorithms are listed in Table \ref{table:2}, which shows that the DCA is faster than the IRL1 in both cases. Therefore, we adopt the former throughout our numerical experiments to solve the problem (\ref{unconstrained}).

\begin{table}[h!]
\centering
 \begin{tabular}{c| c |c } 
 \hline
  &  IRL1 &  DCA \\ [0.5ex] 
 \hline\hline
 Noise Free & 0.169 (0.098) & 0.118 (0.055) \\ 
 \hline
 Noisy & 0.144 (0.071) &  0.124 (0.067) \\ [1ex] 
 \hline
\end{tabular} 
\caption{Computing time of the IRL1 and the DCA algorithms (mean and standard deviation over $10$ outer iterations, when the inner iterations are done for $2N=512$ times).}
\label{table:2}
\end{table}

\subsection{Choice of Parameters}
Next, we carry out a simulation study to understand the effects of the shape parameter $p$ and the scale parameter $\sigma$ on recovering sparse signals based on the generalized error function. In this set of experiments, the true signal is of length $256$ and simulated as $k$-sparse with $k$ in the set $\{2,4,6,\cdots,32\}$. The measurement matrix $A$ is a $64\times 256$ random matrix generated as Gaussian. For each $p$ and $\sigma$, we replicate the experiments $100$ times. It is recorded as one success if the relative error $\lVert \hat{\mathbf{x}}-\mathbf{x}\rVert_2/\lVert \mathbf{x}\rVert_2\leq 10^{-3}$. We show the success rate over the $100$ replicates for various values of parameters $p$ and $\sigma$, while varying the sparsity level $k$.

\begin{figure}
\begin{subfigure}{.5\textwidth}
  \centering
  \includegraphics[width=.9\linewidth]{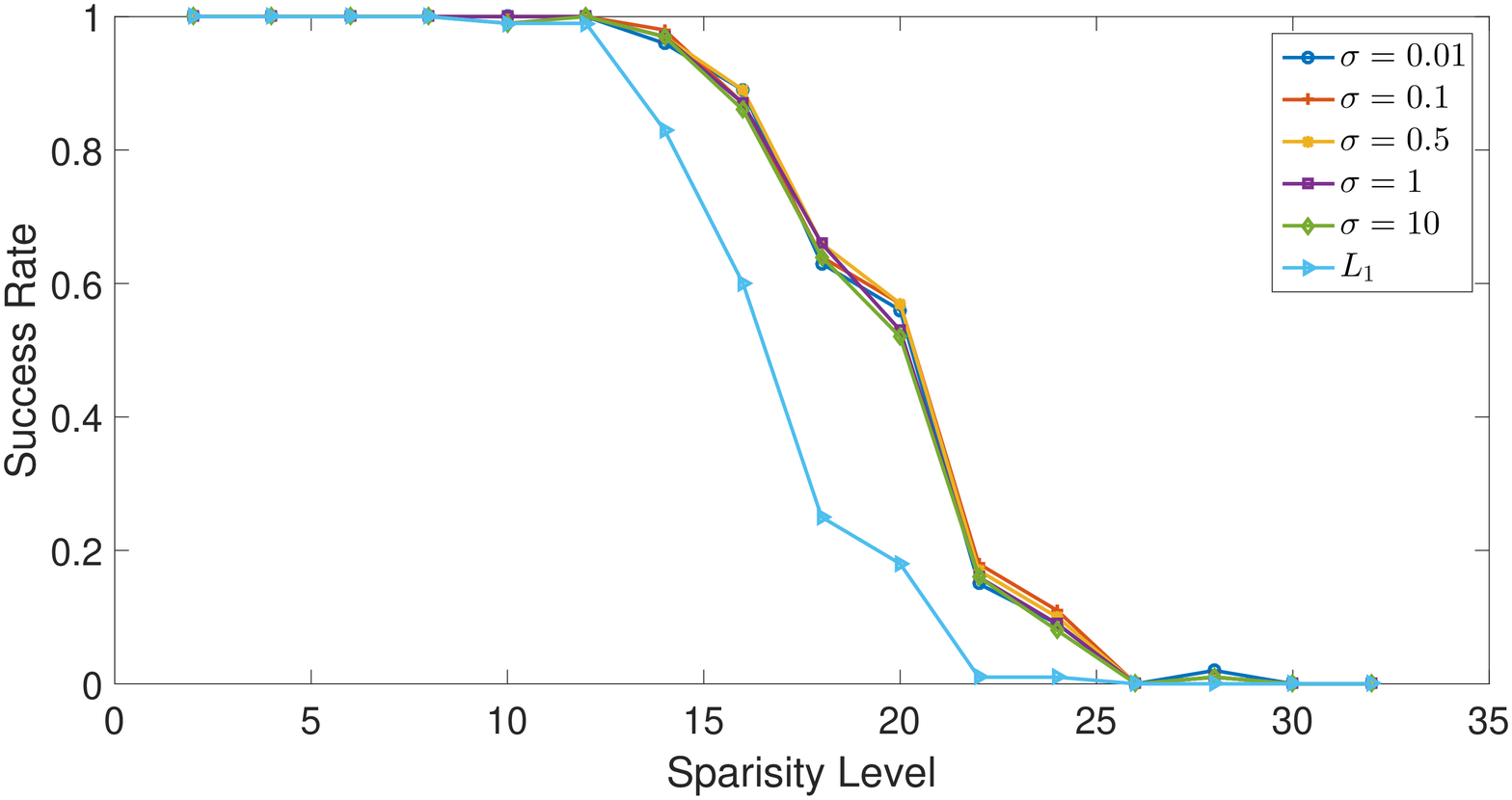}  
  \caption{$p=0.1$}
  \label{fig:sub-first}
\end{subfigure}
\begin{subfigure}{.5\textwidth}
  \centering
  \includegraphics[width=.9\linewidth]{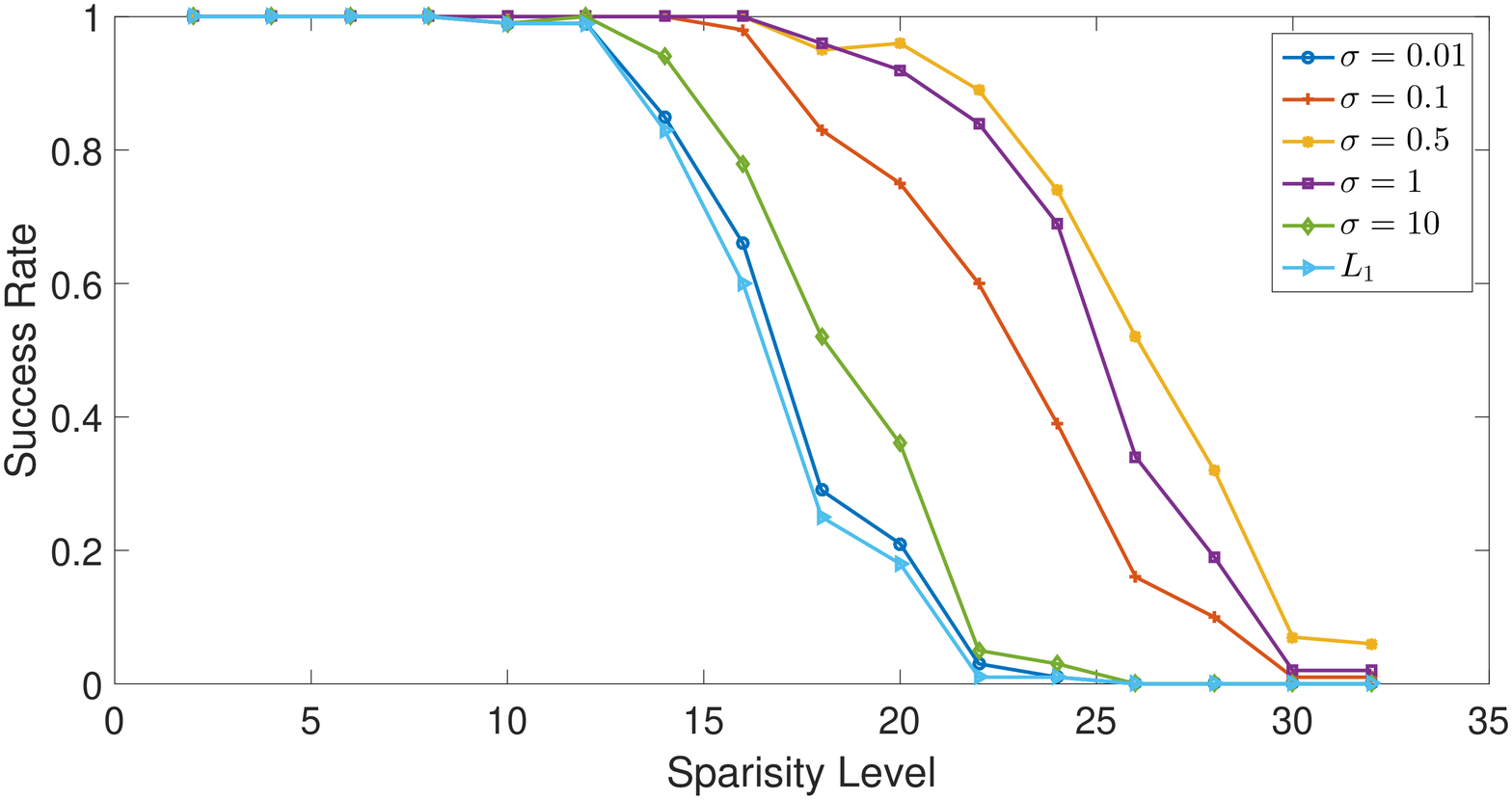}  
  \caption{$p=1$}
  \label{fig:sub-second}
\end{subfigure}

\begin{subfigure}{.5\textwidth}
  \centering
  \includegraphics[width=.9\linewidth]{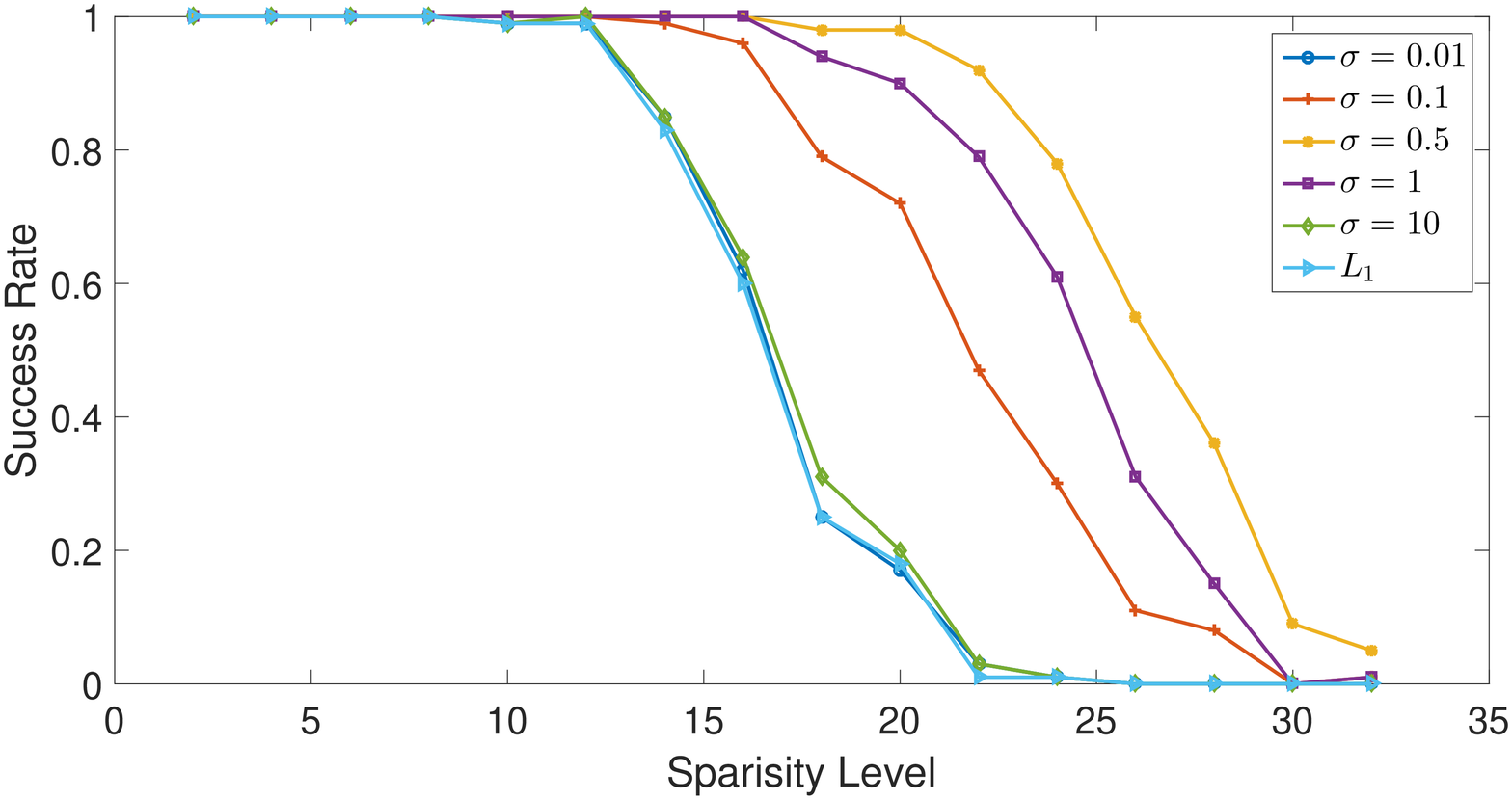}  
  \caption{$p=2$}
  \label{fig:sub-third}
\end{subfigure}
\begin{subfigure}{.5\textwidth}
  \centering
  \includegraphics[width=.9\linewidth]{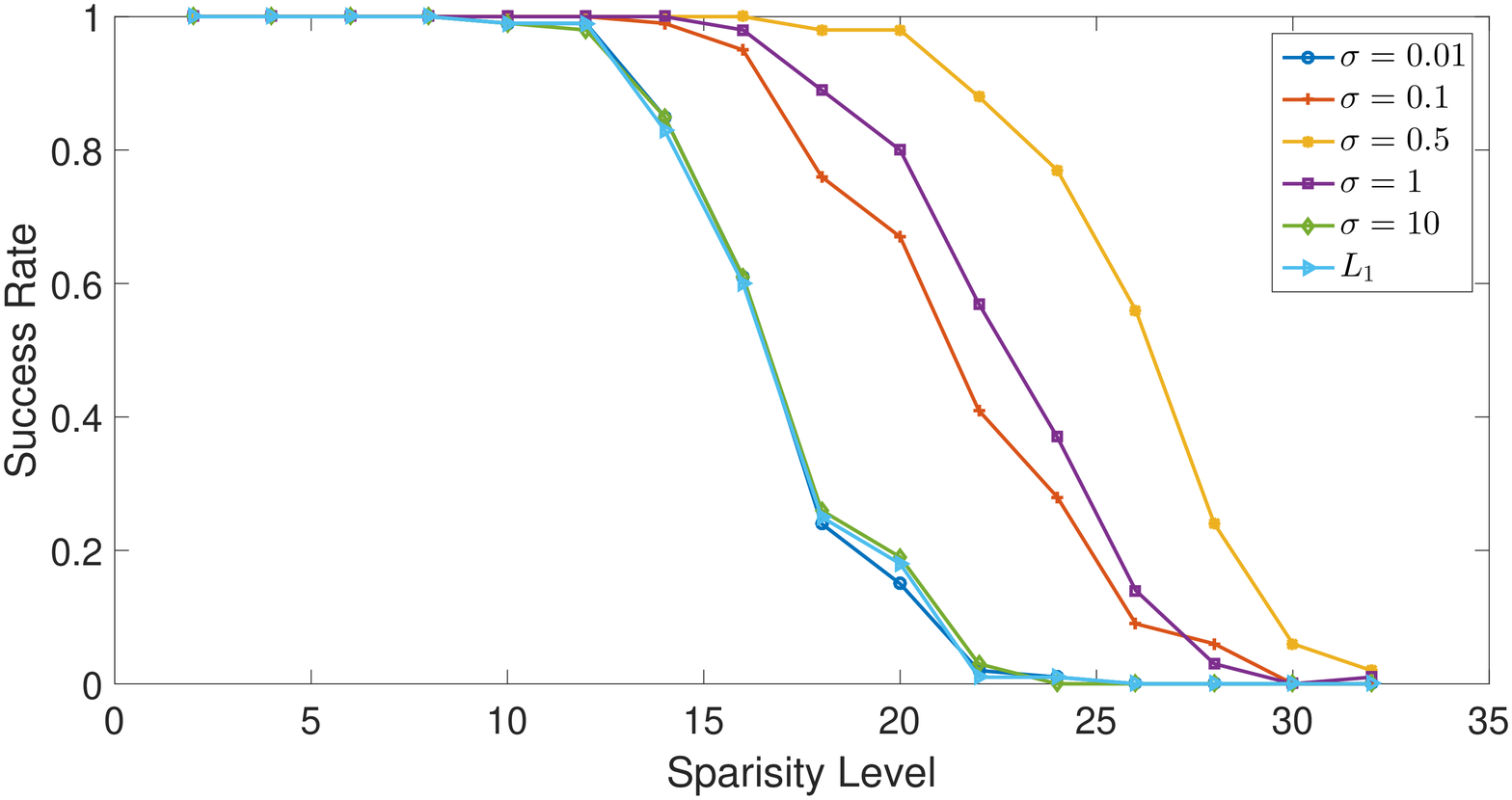}  
  \caption{$p=5$}
  \label{fig:sub-fourth}
\end{subfigure}
\caption{Reconstruction performance of the GERF method with different choices of $\sigma$ for Gaussian random measurements when $p$ is fixed to $0.1,1,2,5$.}
\label{fig:fig}
\end{figure}

As shown in Figure \ref{fig:fig}, the proposed GERF method outperforms the $\ell_1$-minimization (denoted as $L_1$ in the figure) for all tested values of $p$ and $\sigma$. It approaches the $\ell_1$-minimization for a very small value of $\sigma$ (see $\sigma=0.01$ when $p=1,2,5$). In addition, it can be seen from Subfigure \ref{fig:sub-first} that when $p$ is very small (e.g., $p=0.1$), the different choices of $\sigma$ have very limited effect on the recovery performance. Looking at other Subfigures (\ref{fig:sub-second},\ref{fig:sub-third},\ref{fig:sub-fourth}), it is apparent that $\sigma=0.5$ and $\sigma=1$ are the best two choices for a moderate value of $p$ (e.g., $p=1,2,5$).

\begin{figure}
\begin{subfigure}{.5\textwidth}
  \centering
  \includegraphics[width=.9\linewidth]{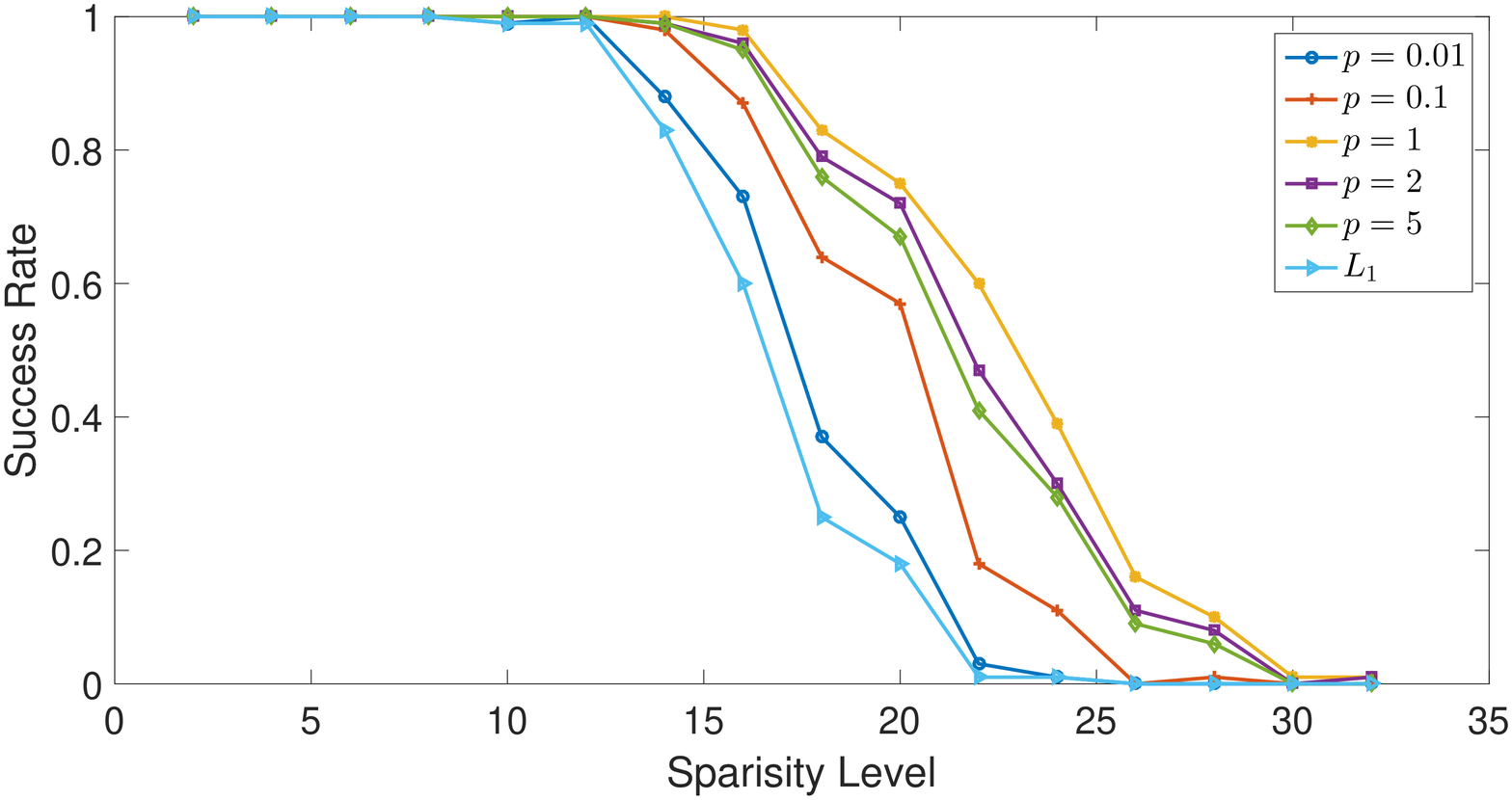}  
  \caption{$\sigma=0.1$}
  \label{fig_sigma:sub-first}
\end{subfigure}
\begin{subfigure}{.5\textwidth}
  \centering
  \includegraphics[width=.9\linewidth]{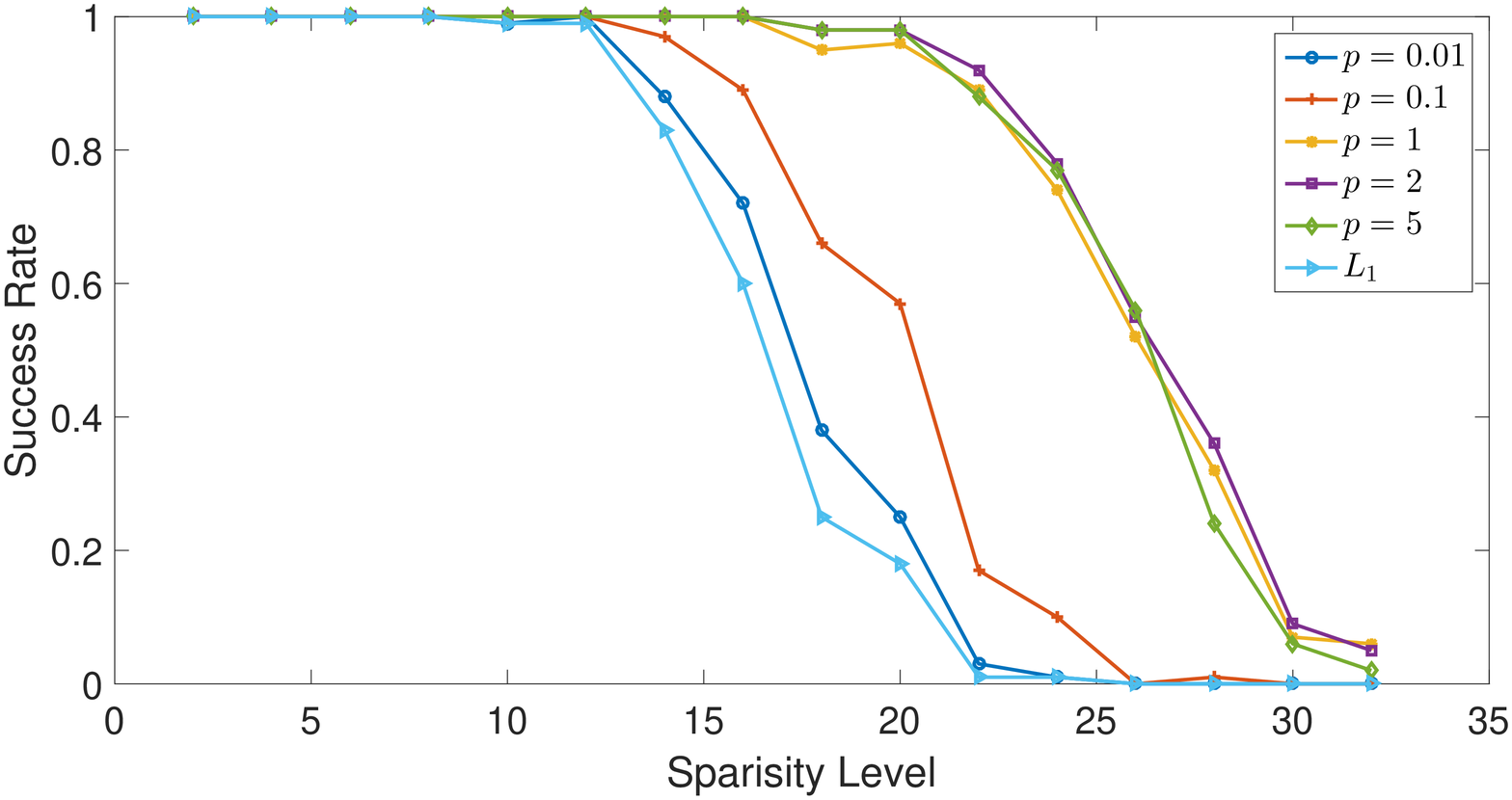}  
  \caption{$\sigma=0.5$}
  \label{fig_sigma:sub-second}
\end{subfigure}

\begin{subfigure}{.5\textwidth}
  \centering
  \includegraphics[width=.9\linewidth]{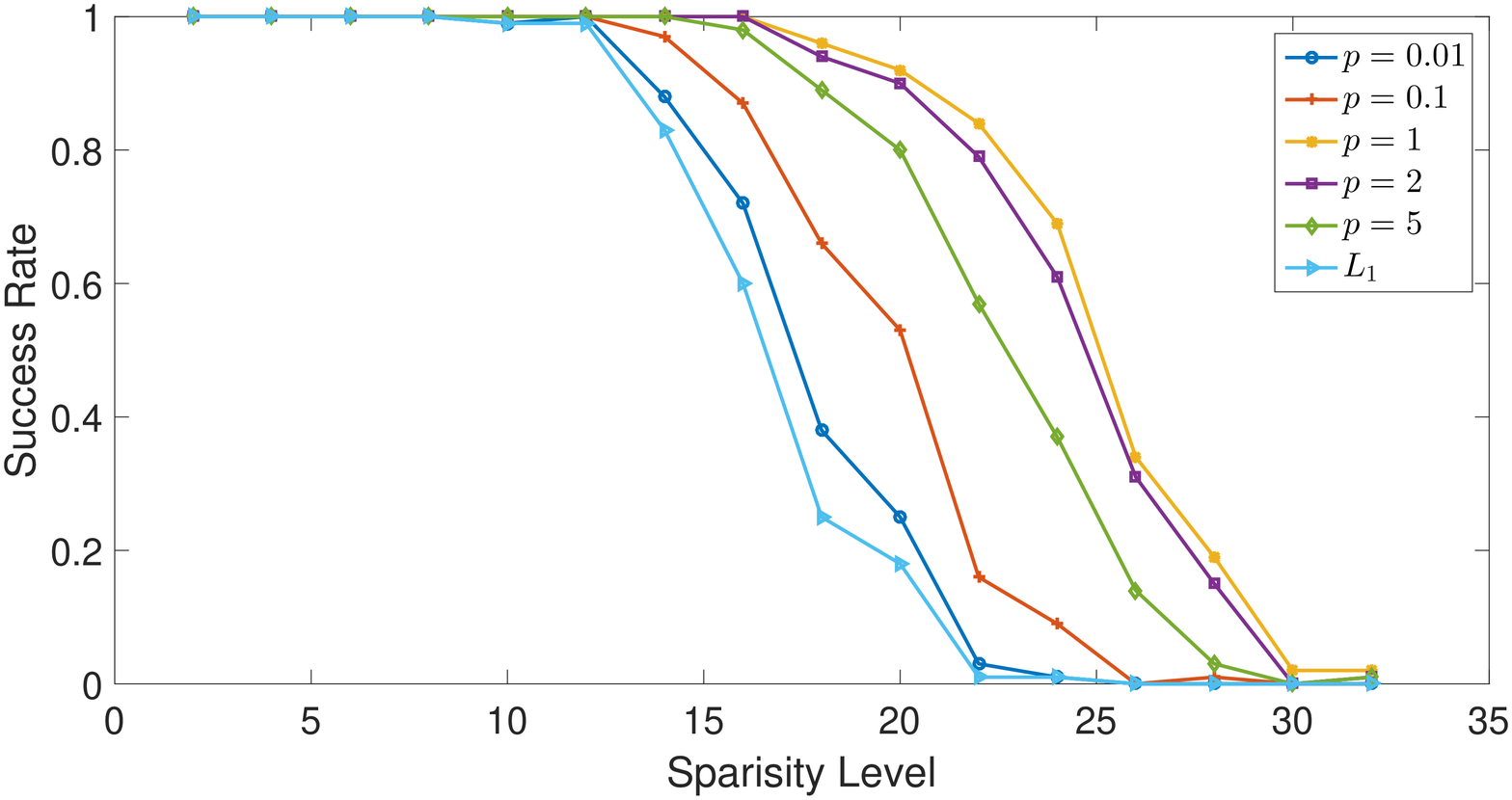}  
  \caption{$\sigma=1$}
  \label{fig_sigma:sub-third}
\end{subfigure}
\begin{subfigure}{.5\textwidth}
  \centering
  \includegraphics[width=.9\linewidth]{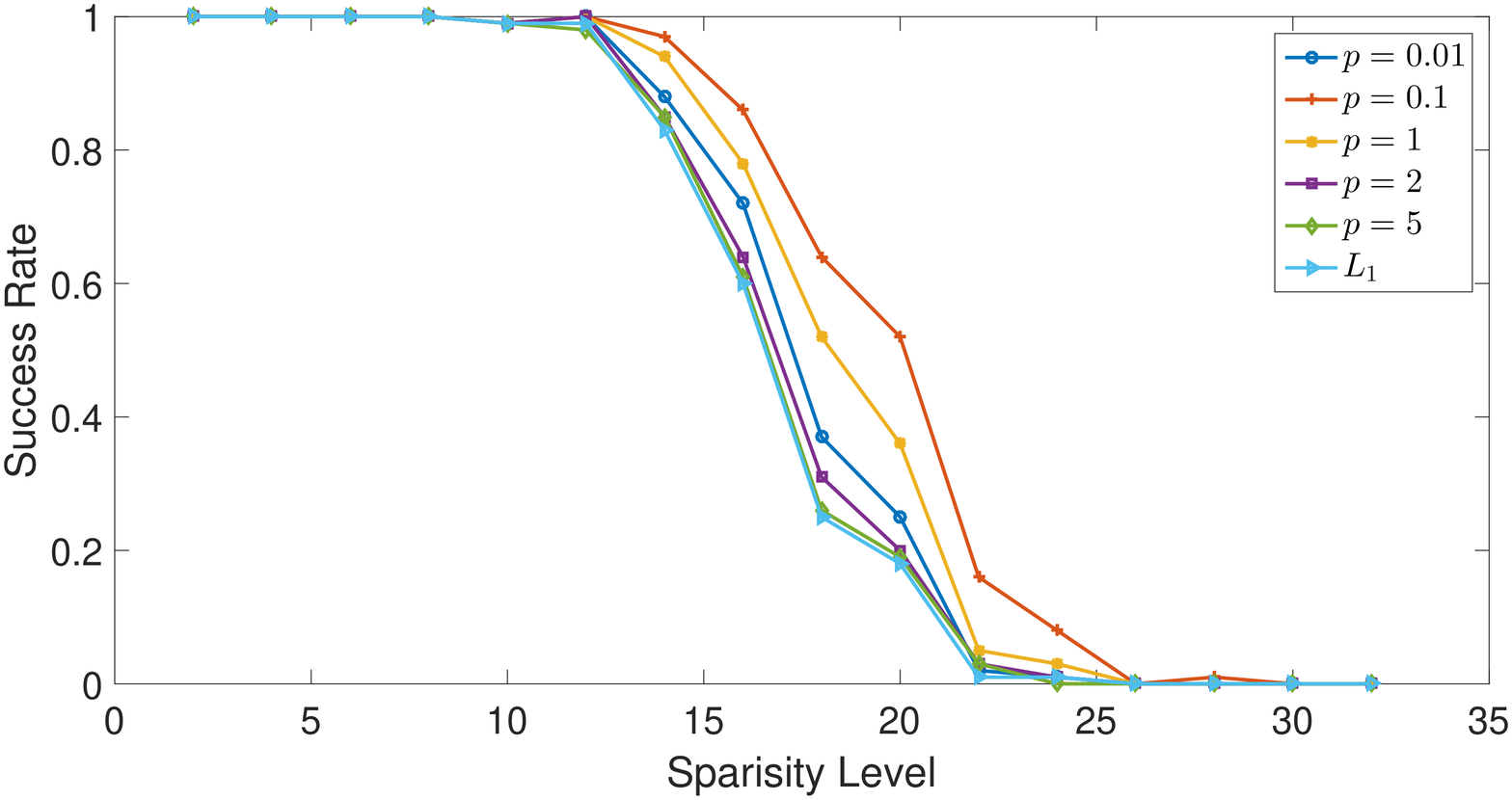}  
  \caption{$\sigma=10$}
  \label{fig_sigma:sub-fourth}
\end{subfigure}
\caption{Reconstruction performance of the GERF method with different choices of $p$ for Gaussian random measurements when $\sigma$ is fixed to $0.1,0.5,1,10$.}
\label{fig_sigma:fig}
\end{figure}

From Figure \ref{fig_sigma:fig}, we see that $p=1$ and $p=2$ are the best two among all tested values of $p$ when $\sigma$ is fixed to $0.1$, $0.5$ and $1$. In contrast, when $\sigma$ is set to be moderate large (e.g., $\sigma=10$), the gap of the recovery performance becomes small for different values of $p$. In this case, $p=0.1$ is the best among all tested values.

\subsection{Algorithm Comparison}
In this subsection, we compare the proposed GERF method with other state-of-the-art sparse recovery methods including ADMM-Lasso, $\ell_{0.5}$, $\ell_1-\ell_2$ and TL1 ($a=1$). In this set of experiments, we consider two scenarios for measurement matrices, i.e., Gaussian random matrix and oversampled discrete cosine transform (DCT) matrix. For the Gaussian random matrix, it is generated as a $64\times 256$ matrix with entries drawn from i.i.d. standard normal distribution. For the oversampled DCT matrix, we use $A=[\mathbf{a}_1,\mathbf{a}_2,\cdots,\mathbf{a}_{1024}]\in\mathbb{R}^{64\times 1024}$ with $\mathbf{a}_j=\frac{1}{\sqrt{64}}\cos\left(\frac{2\pi \mathbf{w}(j-1)}{F}\right), j=1,2,\cdots,1024$, and $\mathbf{w}$ is a random vector uniformly distributed in $[0,1]^{64}$. It is well-known that the DCT matrix is highly coherent and a larger $F$ yields a more coherent matrix.

\begin{figure}[htbp]
	\centering
	\includegraphics[width=\textwidth,height=0.4\textheight,keepaspectratio]{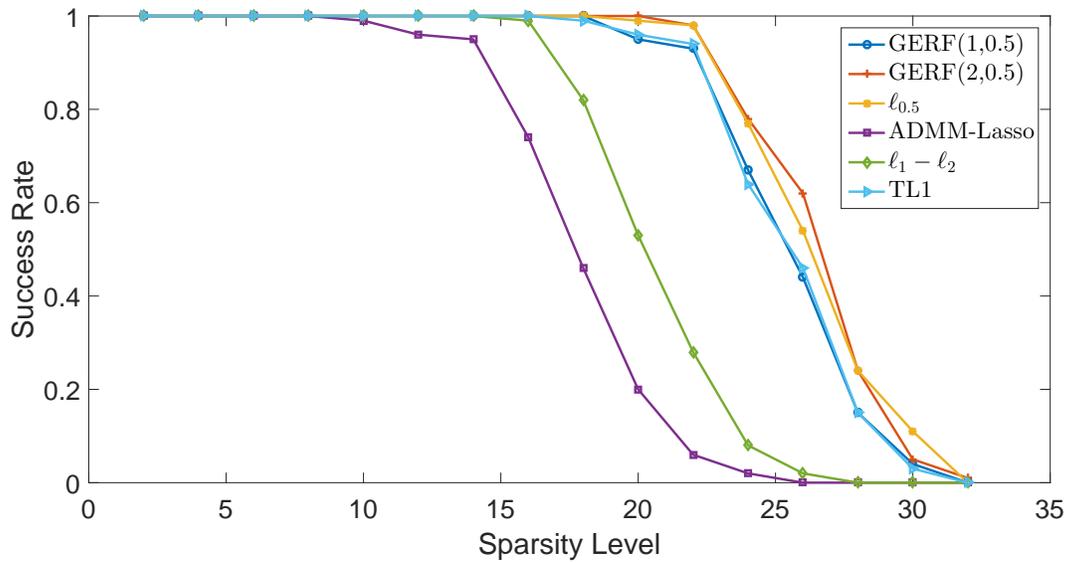}
	\caption{Recovery performance of different algorithms for $64\times 256$ Gaussian random matrices.} \label{algorithm_comparison_gau}
\end{figure}

Figure \ref{algorithm_comparison_gau} shows that for the Gaussian case, the $\ell_{0.5}$ and the GERF with $p=2,\sigma=0.5$ are the best two methods while our proposed GERF method with $p=1,\sigma=0.5$ has an almost same performance as the TL1 method, which are the third best. They are all much better than the ADMM-Lasso and $\ell_1-\ell_2$.

\begin{figure}
\begin{subfigure}{.5\textwidth}
  \centering
  \includegraphics[width=.9\linewidth]{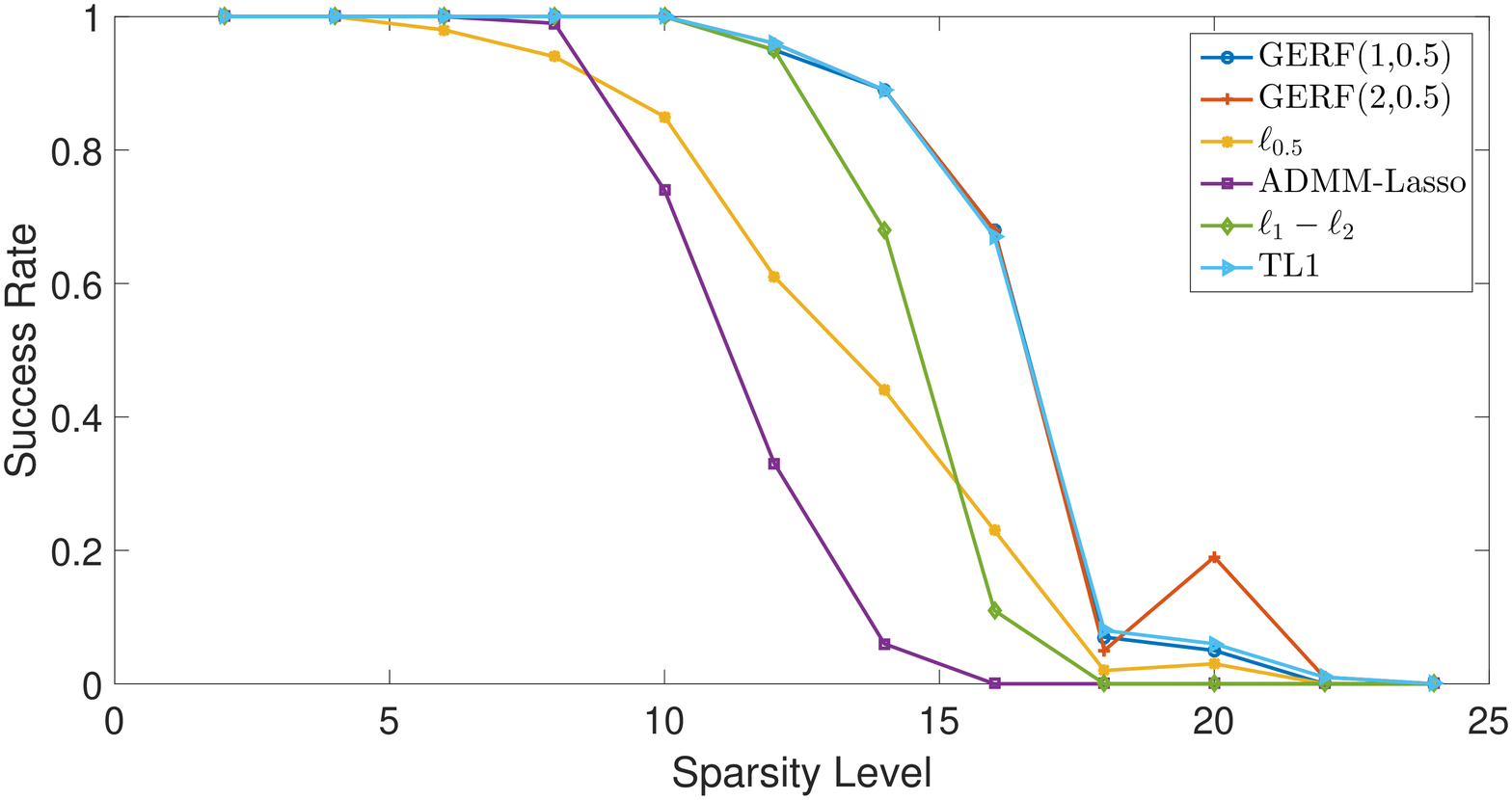}  
  \caption{$F=5$}
  \label{fig_DCT:sub-first}
\end{subfigure}
\begin{subfigure}{.5\textwidth}
  \centering
  \includegraphics[width=.9\linewidth]{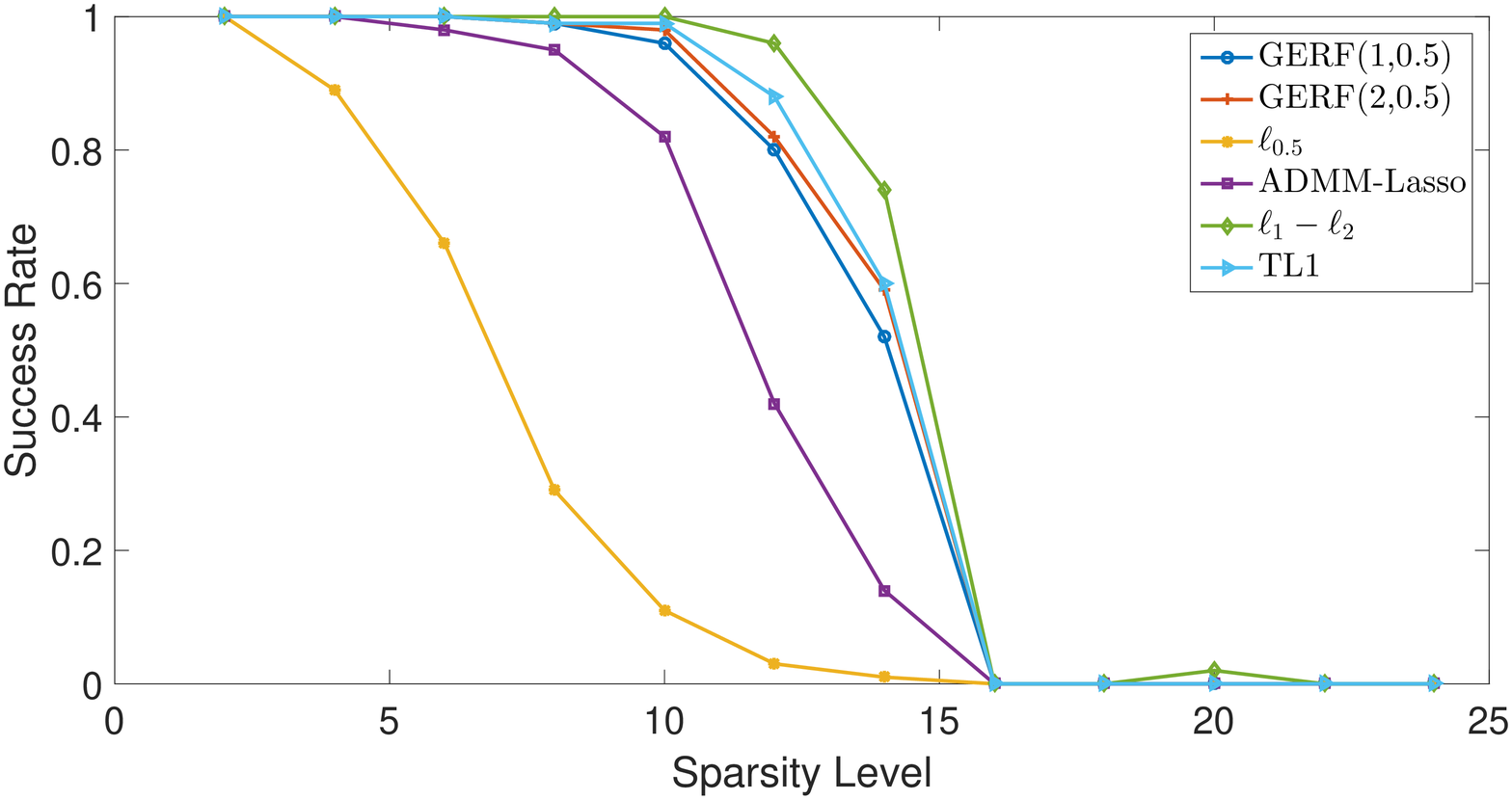}  
  \caption{$F=10$}
  \label{fig_DCT:sub-second}
\end{subfigure}
\caption{Recovery performance of different algorithms for $64\times 1024$ oversampled DCT random matrices. Left panel ($F=5$), Right panel ($F=10$).}
\label{fig_DCT:fig}
\end{figure}

The corresponding results for the oversampled DCT matrix case are displayed in Figure \ref{fig_DCT:fig}. For these two coherent cases, we can observe that although the GERF is not the best, but it is still comparable to the best ones (the TL1 for $F = 5$ and the $\ell_1-\ell_2$ for $F = 10$). Thus, the GERF approach proposed in the present paper gives satisfactory and robust recovery results no matter whether the measurement matrix is coherent or not.

\begin{figure}[htbp]
	\centering
	\includegraphics[width=\textwidth,height=0.4\textheight,keepaspectratio]{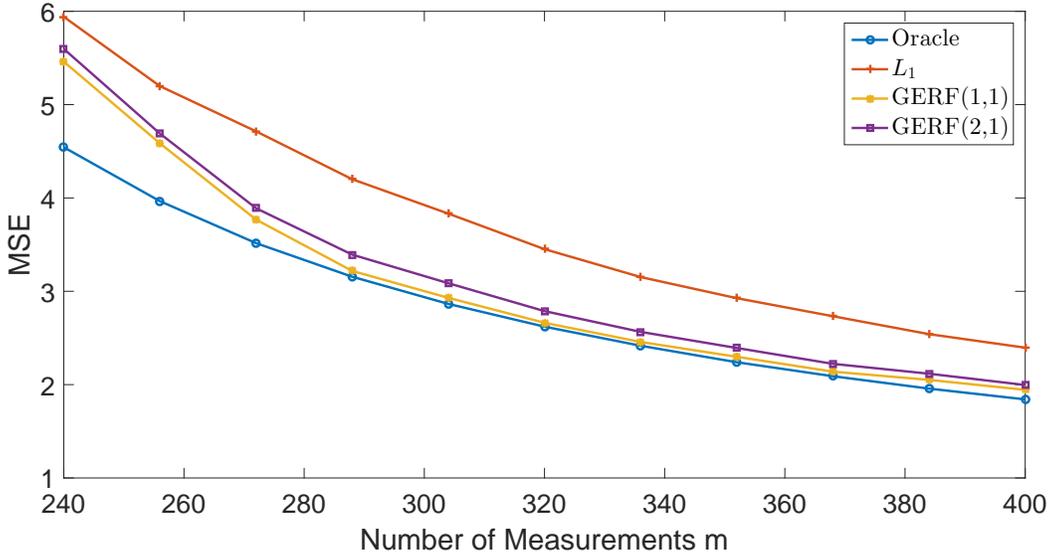}
	\caption{MSE of sparse recovery from noisy Gaussian random measurements while varying the number of measurements $m$.} \label{mse_comparison}
\end{figure}

Furthermore, we present a simulation study for the noisy case. As done in \cite{guo2020novel}, we consider to recover a $130$-sparse signal $\mathbf{x}$ of length $N=512$ from noisy Gaussian random measurements with Gaussian noise (zero mean and standard deviation $\tilde{\sigma}=0.1$). The measurement matrix is of size $m\times 512$ with $m$ ranging from $240$ to $400$. We evaluate the reconstruction performance in terms of the mean-square-error (MSE) over $100$ replications. The benchmark MSE of an oracle solution is computed as $\tilde{\sigma}^2 \mathrm{tr}(A_S^T A_S)^{-1}$, where $S$ is the support of the ground-truth $\mathbf{x}$. In Figure \ref{mse_comparison}, we compare the proposed GERF approaches ($p=1,\sigma=1$ and $p=2,\sigma=1$) with the $L_1$ and the oracle solutions. As we can see, both of our proposed GERF solutions are very close to the oracle solution for large values of $m$, while the $L_1$ solution is obviously biased. 

\subsection{MRI Reconstruction}
\begin{figure}[htbp]
	\centering
	\includegraphics[width=\textwidth,height=0.4\textheight,keepaspectratio]{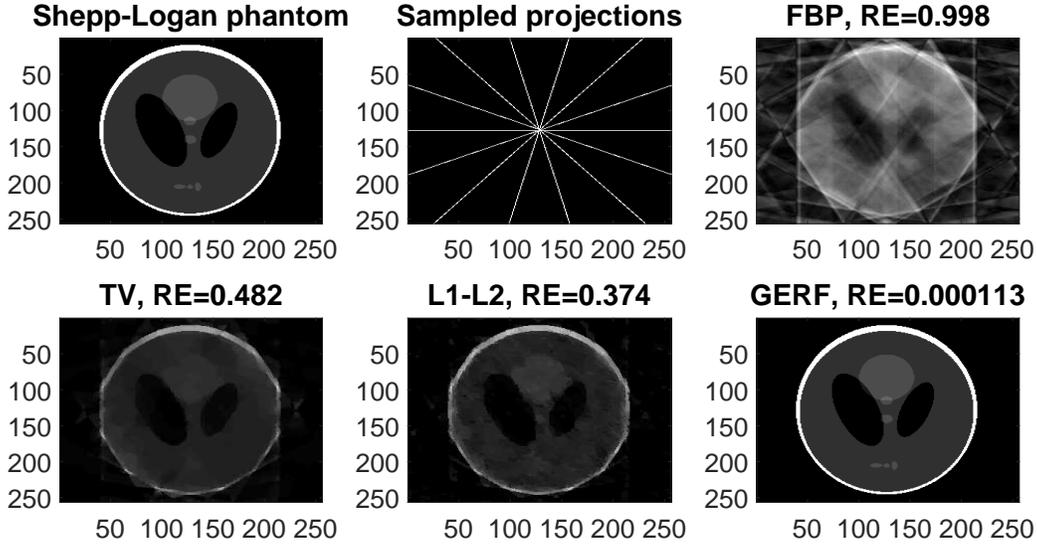}
	\caption{Comparison of RE in MRI reconstruction for the Shepp-Logan phantom of size $256\times 256$ using $7$ projections.} \label{mri_comparison}
\end{figure}

Finally, we examine the particular application of the GERF approach on the gradient in MRI reconstruction. We use the Shepp-Logan phantom of size $256\times 256$ from $7$ radical projections without noise. We compare the reconstruction performance of the filtered back projection (FBP), total variation (TV), L1-L2 \cite{yin2015minimization,lou2015weighted} and our proposed GERF approach with $p=1,\sigma=1$, in terms of the relative error (RE) which equals to $\lVert \hat{u}-u\rVert_{F}/\lVert u\rVert_F$. As we can see from Figure \ref{mri_comparison}, our proposed GERF approach achieves a near-perfect recovery using only $7$ projections, with a relative error $1.13\times 10^{-4}$, which is much better than the classical FBP, TV and L1-L2 methods.

\section{Conclusion}
In this paper, we proposed a concave sparse recovery approach based on the generalized error function, which is less biased than the $L_1$ approach. Theoretically, the recovery analysis was established for both constrained and unconstrained models. Computationally, both IRL1 and DCA algorithms are used to solve the nonsmooth nonconvex problem. A large number of numerical results demonstrated its advantageous performance over the state-of-the-art sparse recovery methods in various scenarios. The MRI phantom image reconstruction test also indicated its superior performance.

In future work, the extensions to structured sparsity (e.g., block sparse recovery and the matrix analogues) are of immense interest. Other practical imaging applications including image denoising and image deblurring should be explored as well. 

\bibliographystyle{unsrt}
\bibliography{references}  


\end{document}